\documentclass{amsart}
\usepackage{amssymb}
\usepackage{graphicx}
\usepackage{xcolor}
\usepackage[hidelinks]{hyperref}
\usepackage[capitalize]{cleveref}
\graphicspath{{figures/}}

\newtheorem{theorem}{Theorem}[section]
\newtheorem{lemma}[theorem]{Lemma}
\theoremstyle{remark}
\newtheorem{remark}[theorem]{Remark}
\numberwithin{equation}{section}
\newcommand{\nn}{\nonumber}
\crefformat{equation}{\textup{#2(#1)#3}}
\Crefname{figure}{Figure}{Figures}
\date{}

\begin{document}

\title[Low-regularity integrators for Schr\"odinger-Poisson]{Low-regularity integrators for the
one-dimensional Schr\"odinger-Poisson equation}

\author{Yue Feng}
\address{School of Mathematics and Statistics, Xi’an Jiaotong University,
Xi'an 710049, China}
\email{yue.feng@xjtu.edu.cn}

\author{Hang Li}
\address{Laboratoire Jacques-Louis Lions, Sorbonne Universit\'e, Paris 75005, France}
\email{hang.li.1@sorbonne-universite.fr}

\author{Wenbo Yang}
\address{School of Mathematics and Statistics, Xi’an Jiaotong University,
Xi'an 710049, China}
\email{wenboyang@stu.xjtu.edu.cn}

\thanks{The authors acknowledge the support from the National Key R$\&$D Program of China (No. 2024YFA1016100), Tianyuan Fund for Mathematics of the National Natural Science Foundation of China (No. 12326102), National Natural Science Foundation of China (No. 12401539) and Jiangsu Provincial Scientific Research Center of Applied Mathematics (No. BK20233002).}

\begin{abstract}
Explicit first- and second-order low-regularity integrators are proposed and analyzed for the one-dimensional Schr\"odinger-Poisson equation with periodic boundary conditions. The schemes are constructed by decomposing the oscillatory terms in Duhamel's formula and exploiting the smoothing effect of the inverse Laplacian in the Poisson potential to reduce the regularity requirements. For $r>1/2$, the first-order method converges in $H^r$ for solutions bounded in $H^r$, without requiring any additional spatial regularity. The second-order method converges in $H^r$ for solutions bounded in $H^{r+1}$. In addition, first- and second-order convergence in $L^2$ is established for solutions bounded in $H^{1/2}$ and $H^{5/4}$, respectively. Compared with classical splitting methods and exponential integrators, the proposed schemes achieve the corresponding convergence orders under weaker regularity assumptions. Numerical experiments confirm the theoretical error estimates and demonstrate the favorable performance of the proposed schemes for low-regularity initial data.
\end{abstract}

\maketitle

\section{Introduction}
The Schr\"odinger--Poisson equation (SPE) arises in a broad range of applications, including quantum mechanics \cite{BEGMY,LB}, nonlinear optics \cite{BSMRS,POM}, cosmology \cite{SNE,VN}, and semiconductor physics and nanoelectronics \cite{MRS,TGPR,WSD}. In particular, it provides a model for electron transport in quantum wells, superlattices and nanowires.
In this paper, we consider the one-dimensional SPE on the torus $\mathbb{T}=\mathbb{R}/(2\pi\mathbb{Z})$ with periodic boundary conditions:
\begin{equation}
\label{eq:1D}
\begin{cases}
i\partial_t\psi(x,t)
=-\frac{1}{2}\Delta\psi(x,t)+V(x,t)\psi(x,t),
& x\in\mathbb{T},\quad t>0,\\
-\Delta V(x,t)=|\psi(x,t)|^2-C,
& x\in\mathbb{T},\quad t\geq0,\\
\psi(x,0)=\psi_0(x),
& x\in\mathbb{T}.
\end{cases}
\end{equation}
Here, $t$ and $x$ denote the temporal and spatial
variables, respectively, 
and $\psi(x,t)\in\mathbb{C}$ is the wave function
with prescribed initial data $\psi_0(x)$.
The real-valued electric potential $V(x,t)$ is
determined self-consistently by the Poisson equation. The normalization
$\int_{\mathbb{T}}V(x,t)\,dx=0$ guarantees uniqueness of the potential, while the constant $C=\frac{1}{2\pi}\int_{\mathbb{T}}|\psi_0(x)|^2\,dx$
ensures compatibility with the Poisson equation.

The SPE \cref{eq:1D} possesses several fundamental
dynamical properties \cite{ABB,BW}. In particular, it is invariant under the combined transformation of time reversal and complex conjugation.
Moreover, the mass \cite{AKKM,ZY}
\begin{equation}
N(t):=\int_{\mathbb{T}}|\psi(x,t)|^2\,dx
\equiv
\int_{\mathbb{T}}|\psi_0(x)|^2\,dx
=N(0),\qquad t\geq0,
\end{equation}
and the energy
\begin{align}
E(t)
&:=\int_{\mathbb{T}}
\left[
\frac{1}{2}|\partial_x\psi(x,t)|^2
+\frac{1}{2}V(x,t)|\psi(x,t)|^2
\right]\,dx \nn\\
&\equiv\int_{\mathbb{T}}
\left[
\frac{1}{2}|\partial_x\psi_0(x)|^2
+\frac{1}{2}V(x,0)|\psi_0(x)|^2
\right]\,dx
=E(0),\qquad t\geq0,
\end{align}
are conserved.

The SPE has been studied extensively from both analytical and numerical perspectives. Analytical results concerning existence, uniqueness and asymptotic behavior of solutions can be found in \cite{BM,Cas,IZ,ZZ} and the references therein. From the numerical viewpoint, a variety of accurate and efficient methods have been developed, including finite difference
methods \cite{AKKM,RS,ZY}, time-splitting methods
\cite{AKKT,BMS,LC,SLZMW}, exponential time differencing methods
\cite{Lu}, and Gaussian beam methods \cite{JWY}. Nevertheless, standard time discretizations generally require relatively high spatial 
regularity to attain their optimal convergence orders and may suffer from order reduction for low-regularity initial data.

In recent years, low-regularity integrators have been developed to address this difficulty for nonlinear Schr\"odinger equations and related dispersive problems. Ostermann and Schratz \cite{OS} introduced low-regularity exponential-type integrators for semilinear Schr\"odinger equations. 
The first-order method for the cubic nonlinear Schr\"odinger equation requires only one additional spatial derivative because the dominant oscillatory terms are integrated exactly. Higher-order methods
were subsequently developed in \cite{BS,KOS,OWY}, and integrators under
weaker regularity assumptions were investigated in \cite{CLL,ORS,WY}.
Fully discrete error analyses were established in
\cite{LW,OY}, and extensions to nonperiodic
boundary conditions were considered in \cite{AB,BLW}. Low-regularity integrators have also been constructed for other
dispersive equations, including the KdV equation \cite{HS,WZKdV},
nonlinear Klein--Gordon equation \cite{WZKG}, and nonlinear
Dirac equation \cite{SWZ}.

To the best of our knowledge, low-regularity integrators for the one-dimensional periodic SPE have not yet been studied and analyzed. In this paper, we construct and analyze explicit first- and second-order
low-regularity integrators for this problem.
A direct adaptation of the dominant-phase approximation from \cite{OS}
leads to a first-order remainder estimate that requires one additional
spatial derivative. Consequently, it does not improve the regularity requirement of the Lie--Trotter splitting method. To avoid this loss of regularity, we exploit the smoothing effect of the inverse Laplacian in the Poisson potential, decompose the oscillatory phase in a suitable way, and select a different factor for exact integration. This choice allows the spatial derivative appearing in the remainder to be compensated by the inverse Laplacian. For the resulting first-order integrator, we prove first-order convergence in $H^r$ for solutions bounded in $H^r$, with $r>1/2$, and therefore without any additional spatial regularity.
We then extend the construction to obtain a second-order
method and prove second-order convergence in $H^r$ for solutions bounded in $H^{r+1}$. For the $L^2$-norm, we establish first-order convergence for solutions bounded in $H^{1/2}$ and second-order convergence for solutions bounded in $H^{5/4}$.
Numerical experiments support the theoretical error estimates and
illustrate the advantages of the proposed schemes for low-regularity
initial data.

The remainder of this paper is organized as follows.
In Section~2, we construct the first-order low-regularity
integrator and establish its stability and global error bounds in $H^r$ ($r>1/2$)
and $L^2$. In Section~3, we extend the construction to
a second-order method and prove the corresponding convergence results.
Numerical experiments are presented in Section~4
to illustrate the theoretical results and compare
the proposed schemes with classical numerical methods. Finally, some conclusions are drawn in Section~5.

\section{A first-order low-regularity integrator}
In this section, we derive a first-order low-regularity integrator
for the SPE \eqref{eq:1D} and establish its error estimate.
We first introduce some notation used throughout the paper. For $f, g \in L^2(\mathbb{T})$, we denote the Fourier expansion of $f$ by
\[
f(x)=\sum_{k\in\mathbb{Z}}\widehat f_k e^{ikx},
\quad
\widehat f_k =\frac1{2\pi}\int_{\mathbb{T}} f(x)e^{-ikx}dx,
\]
and the normalized $L^2$ inner product by
\[
 \langle f,g\rangle
 := \frac{1}{2\pi}\int_{\mathbb T}f(x)\overline{g(x)}\,dx.
\]
We write $\langle k\rangle=(1+|k|^2)^{1/2}$ and equip
the Sobolev space $H^s(\mathbb T)$ with the norm
\[
\|f\|_{H^s}^2
:=\sum_{k\in\mathbb Z}\langle k\rangle^{2s}|\widehat f_k|^2.
\]
For any integer $m\geq1$ and any $\alpha\in\mathbb R$,
we define the Fourier multipliers $\partial_x^{-m}$
and $|\partial_x|^\alpha$ by
\begin{equation*}
\widehat{\partial_x^{-m}f}_k=
\begin{cases} 
(ik)^{-m} \widehat{f}_k,&k \neq 0, \\
0, & k = 0,
\end{cases}
\qquad
\widehat{|\partial_x|^\alpha f}_k=
\begin{cases}
|k|^\alpha\widehat f_k, & k\neq0,\\
0, & k=0.
\end{cases}
\end{equation*}
In particular, we have
$(-\Delta)^{-1}=-\partial_x^{-2}$ in one dimension.

Throughout this paper, $c$ denotes a generic positive constant independent of $\tau$ and $n$, whose value may vary from one occurrence to another. We write $c(M)$ when its dependence on $M$ needs to be emphasized. The notation $A\lesssim B$ means that $A\le cB$, and $O_{H^r}(s^p)$ denotes a remainder
whose $H^r$-norm is bounded by $cs^p$.

{Assuming that the exact solution exists on $[0,T]$, we write}
\begin{equation}\label{eq:Mr}
 M_r=\sup_{0\leq t\leq T}\|\psi(t)\|_{H^r}.
\end{equation}
\subsection{Construction of the scheme}
Choose a time step size $\tau > 0$ and define the time steps as $t_n=n\tau$ for $n=0,1,\dots$. In this paper, we shall abbreviate $\psi(\cdot, t)$ by $\psi(t)$ whenever no confusion can arise. By Duhamel's principle, the mild solution $\psi(t)$ of the SPE \eqref{eq:1D} is 
\begin{equation}
\label{eq:Duhamel1}
\psi(t_n+\tau)=e^{i\tau\frac{\Delta}{2}}\psi(t_n)-i\int_{0}^{\tau} e^{i(\tau-s)\frac{\Delta}{2}}V(t_n+s)\psi(t_n+s)ds.
\end{equation}

For the estimates below, let $r>1/2$ and assume that
$\psi\in L^\infty([0,T];H^r(\mathbb T))$. Since \(e^{is\frac{\Delta}{2}}\) is an isometry on \(H^r\) and \(H^r\) is an algebra for \(r>1/2\), we have
\begin{equation*}
\|\psi(t_n+s)-e^{is\frac{\Delta}{2}}\psi(t_n)\|_{H^{r}}\leq \int_0^s \left\|V(t_n+\sigma)\psi(t_n+\sigma)\right\|_{H^r}d\sigma\leq c(M_r)s.
\end{equation*}
This estimate allows us to approximate the solution within one time step by
\begin{equation}
\label{eq:approx1}
\psi(t_n+s)=e^{is\frac{\Delta}{2}}\psi(t_n)+O_{H^r}(s), \quad 0\le s\le\tau.
\end{equation}
The approximation \eqref{eq:approx1} yields the following estimate
\[
\begin{aligned}
&\left\|
V(t_n+s)
-(-\Delta)^{-1}
\left(\left|e^{is\frac{\Delta}{2}}\psi(t_n)\right|^2\right)
\right\|_{H^r}\\
&\lesssim
\left(
\|\psi(t_n+s)\|_{H^r}
+\left\|e^{is\frac{\Delta}{2}}\psi(t_n)\right\|_{H^r}
\right)
\left\|\psi(t_n+s)-e^{is\frac{\Delta}{2}}\psi(t_n)\right\|_{H^r}\\
&\le c(M_r)s.
\end{aligned}
\]
Combining the estimates above, we have
\[
\begin{aligned}
&\left\|
V(t_n+s)\psi(t_n+s)
-(-\Delta)^{-1}
\left(\left|e^{is\frac{\Delta}{2}}\psi(t_n)\right|^2\right)
e^{is\frac{\Delta}{2}}\psi(t_n)
\right\|_{H^r}\\
&\lesssim
\|V(t_n+s)\|_{H^r}
\left\|\psi(t_n+s)-e^{is\frac{\Delta}{2}}\psi(t_n)\right\|_{H^r}\\
&\quad+
\left\|
V(t_n+s)
-(-\Delta)^{-1}
\left(\left|e^{is\frac{\Delta}{2}}\psi(t_n)\right|^2\right)
\right\|_{H^r}
\left\|e^{is\frac{\Delta}{2}}\psi(t_n)\right\|_{H^r}\\
&\le c(M_r)s.
\end{aligned}
\]
Therefore, substituting \eqref{eq:approx1} into \eqref{eq:Duhamel1}, we obtain
\begin{equation}
\label{eq:expansion1}
\psi(t_n+\tau)
=e^{i\tau\frac{\Delta}{2}}\psi(t_n)
-i e^{i\tau\frac{\Delta}{2}}I^\tau(\psi(t_n))
+O_{H^r}(\tau^2),
\end{equation}
where, for a given function $\psi$, we define
\begin{equation}
\label{eq:I^tau}
I^{\tau}(\psi):=\int_{0}^{\tau} e^{-is\frac{\Delta}{2}} (-\Delta)^{-1}\left(\left|e^{is\frac{\Delta}{2}}\psi\right|^2\right)e^{is\frac{\Delta}{2}}\psi \,ds.
\end{equation}

Expressing \eqref{eq:I^tau} in Fourier space, we obtain
\begin{equation}
\label{eq:I^tau_expansion}
I^\tau(\psi)=\sum_{k\in \mathbb{Z}}\sum_{\substack{k_1,k_2,k_3\in \mathbb{Z}\\k=-k_1+k_2+k_3\\k_1\neq k_2}}\frac{1}{(k_2-k_1)^2}\int_{0}^{\tau}e^{i\frac{s}{2}(k^2+k_1^2-k_2^2-k_3^2)}ds \overline{\widehat{\psi}_{k_1}}\widehat{\psi}_{k_2}\widehat{\psi}_{k_3}e^{ikx}.
\end{equation}

Using \(k=-k_1+k_2+k_3\), we rewrite the integrand as
\begin{equation}
\label{eq:integrand_decomposition1}
\begin{aligned}
e^{i\frac{s}{2}(k^2+k_1^2-k_2^2-k_3^2)}
=e^{isk_3(k_2-k_1)}+e^{isk_3(k_2-k_1)}
\bigl(e^{isk_1(k_1-k_2)}-1\bigr).
\end{aligned}
\end{equation}
The first term yields an approximation that can be efficiently
evaluated in physical space, while the second contributes a
remainder that can be controlled without additional spatial
regularity. Accordingly, we write
\begin{equation}
\label{eq:I^tau_decomposition1}
I^\tau(\psi)=I_1^\tau(\psi)+R_1^\tau(\psi),
\end{equation}
where
\begin{equation}
\label{eq:I1}
I_1^\tau(\psi):=
\sum_{k\in\mathbb Z}
\sum_{\substack{k_1,k_2,k_3\in\mathbb Z\\
k=-k_1+k_2+k_3\\k_1\ne k_2}}
\frac{\overline{\widehat\psi_{k_1}}
\widehat\psi_{k_2}\widehat\psi_{k_3}}
{(k_2-k_1)^2}
\int_0^\tau e^{isk_3(k_2-k_1)}\,dse^{ikx},
\end{equation}
and
\begin{align}
\label{eq:remainder1}
R_1^\tau(\psi)
=\sum_{k\in\mathbb Z}
\sum_{\substack{k_1,k_2,k_3\in\mathbb Z\\
k=-k_1+k_2+k_3\\ k_1\neq k_2}}
\frac{\overline{\widehat\psi_{k_1}}
\widehat\psi_{k_2}\widehat\psi_{k_3}}
{(k_2-k_1)^2}
\int_0^\tau e^{is k_3(k_2-k_1)}
\left(e^{is k_1(k_1-k_2)}-1\right)ds \,e^{ikx}.
\end{align}
To evaluate $I_1^\tau(\psi)$, we separate the contributions
from $k_3=0$ and $k_3\ne0$, denoted by
$I_{1,1}^\tau(\psi)$ and $I_{1,2}^\tau(\psi)$,
respectively. Thus,
\[
I_1^\tau(\psi)=I_{1,1}^\tau(\psi)+I_{1,2}^\tau(\psi).
\]
For $k_3=0$, the corresponding contribution reduces to
\begin{align}
\label{eq:I^tau11}
I^{\tau}_{1,1}(\psi)=\sum_{k\in \mathbb{Z}}\sum_{\substack{k_1,k_2\in \mathbb{Z}\\k=-k_1+k_2\\k_1\neq k_2}}\tau\frac{1}{(k_2-k_1)^2}\overline{\widehat{\psi}_{k_1}}\widehat{\psi}_{k_2}\widehat{\psi}_{0}e^{ikx}=-\widehat{\psi}_0\tau\partial_x^{-2}(\overline{\psi}\psi).
\end{align}
For $k_3\neq 0$, by using the relation $2k_3(k_2-k_1)=k^2-k_3^2-(k_2-k_1)^2$, the integral becomes
\begin{equation*}
\int_{0}^{\tau}e^{i\frac{s}{2} \left[2k_3(k_2-k_1)\right]}ds=\frac{e^{i\frac{\tau}{2}(k^2-k_3^2-(k_2-k_1)^2)}-1}{ik_3(k_2-k_1)}.
\end{equation*}
Using this expression, $I^{\tau}_{1,2}(\psi)$ can be written in the physical space as
\begin{align}
\label{eq:I^tau12}
I^{\tau}_{1,2}(\psi)=-ie^{-i\frac{\tau}{2}\Delta}\left[\left(e^{i\frac{\tau}{2}\Delta}\partial_x^{-3}(\overline{\psi}\psi)\right)\left(e^{i\frac{\tau}{2}\Delta}\partial_x^{-1}\psi\right)\right]+i\partial_x^{-3}(\overline{\psi}\psi)\partial_x^{-1}\psi.
\end{align}

Let $\psi^n$ be the numerical approximation of $\psi(t_n)$ for $n \geq0$. Combining \eqref{eq:expansion1}, \eqref{eq:I^tau_decomposition1}, \eqref{eq:I^tau11} and \eqref{eq:I^tau12}, and dropping the remainder terms, we obtain the first-order low-regularity integrator (LRI1) for the SPE \eqref{eq:1D}
\begin{align}
\label{eq:LRI1}
\psi^{n+1} &= e^{i\frac{\tau}{2}\Delta}\psi^n-\left(e^{i\frac{\tau}{2}\Delta}\partial_x^{-3}(\overline{\psi^n}\psi^n)\right)\left(e^{i\frac{\tau}{2}\Delta}\partial_x^{-1}\psi^n\right)\\
&\quad+e^{i\frac{\tau}{2}\Delta}\left(\partial_x^{-3}(\overline{\psi^n}\psi^n)\partial_x^{-1}\psi^n\right)+i\tau\widehat{\psi^n}_0e^{i\frac{\tau}{2}\Delta}\partial_x^{-2}(\overline{\psi^n}\psi^n),\nonumber
\end{align}
where $\psi^0=\psi_0$ and $\widehat{\psi^n}_0$ denotes the zero Fourier mode of $\psi^n$.
We denote the numerical flow in \eqref{eq:LRI1} by $\Phi_\tau^{\mathrm{LRI1}}$, i.e.,
$\psi^{n+1}=\Phi_\tau^{\mathrm{LRI1}}(\psi^n)$. The resulting LRI1 scheme is fully explicit. 

\begin{remark}
Following the general construction in \cite{FMW}, we obtain an explicit symmetric two-step method (sLRI1)
\begin{equation*}
 \psi^{n+1}
 =e^{i\tau\Delta}\psi^{n-1}
 +\Phi_\tau^{\mathrm{LRI1}}(\psi^n)
 -e^{i\tau\Delta}\Phi_{-\tau}^{\mathrm{LRI1}}(\psi^n),
 \quad n\ge1.
\end{equation*}
Here, $\Phi_{-\tau}^{\mathrm{LRI1}}$ is obtained from \eqref{eq:LRI1}
by replacing $\tau$ with $-\tau$.
Substituting \eqref{eq:LRI1} into the above formula, we have
\begin{equation}
\label{eq:sLRI1}
\begin{aligned}
 \psi^{n+1}
 &=e^{i\tau\Delta}\psi^{n-1}-\left(
 e^{i\tau\frac{\Delta}{2}}\partial_x^{-3}(|\psi^n|^2)
 \right)
 \left(
 e^{i\tau\frac{\Delta}{2}}\partial_x^{-1}\psi^n
 \right)\\
 &\quad+e^{i\tau\Delta}\left[
 \left(
 e^{-i\tau\frac{\Delta}{2}}\partial_x^{-3}(|\psi^n|^2)
 \right)
 \left(
 e^{-i\tau\frac{\Delta}{2}}\partial_x^{-1}\psi^n
 \right)
 \right]\\
 &\quad +2i\tau\widehat{\psi^n}_0
 e^{i\tau\frac{\Delta}{2}}\partial_x^{-2}(|\psi^n|^2),
 \quad n\ge1,
\end{aligned}
\end{equation}
with
\[
 \psi^1=\Phi_\tau^{\mathrm{LRI1}}(\psi^0),
 \quad \psi^0=\psi_0.
\]
{The sLRI1 scheme is included for numerical comparison,
while the convergence analysis in this paper concerns LRI1 and LRI2.}
\end{remark}

\subsection{Convergence analysis}
In this subsection, we establish stability and local error estimates
for the LRI1 scheme \eqref{eq:LRI1} and then derive the global error bound.
\begin{lemma}[Stability]
\label{lem:stability1}
Let $r>1/2$ and $f,g\in H^r(\mathbb{T})$. We have 
\[
\|\Phi_\tau^{\mathrm{LRI1}}(f)-\Phi_\tau^{\mathrm{LRI1}}(g)\|_{H^r}
\leq \left(1+\tau c(\|f\|_{H^r},\|g\|_{H^r})\right)\|f-g\|_{H^r}.
\]
\end{lemma}
\begin{proof}
By the Fourier representation of $I_1^\tau$, we have
\[
\widehat{I_1^\tau(f)}_k
=
\sum_{\substack{k_1,k_2,k_3\in \mathbb{Z}\\k=-k_1+k_2+k_3\\k_1\neq k_2}}
\frac{\overline{\widehat f_{k_1}}\widehat f_{k_2}\widehat f_{k_3}}
     {(k_2-k_1)^2}
\int_0^\tau e^{is k_3(k_2-k_1)}\,ds.
\]
Noting that
\[
\left|
\int_0^\tau e^{is k_3(k_2-k_1)}\,ds
\right|
\le\tau,
\]
and using the identity
\[
\begin{aligned}
&\overline{\widehat f_{k_1}}\widehat f_{k_2}\widehat f_{k_3}
-\overline{\widehat g_{k_1}}\widehat g_{k_2}\widehat g_{k_3}\\
&=
\overline{\widehat{f-g}_{k_1}}\widehat f_{k_2}\widehat f_{k_3}
+\overline{\widehat g_{k_1}}\widehat{f-g}_{k_2}\widehat f_{k_3}
+\overline{\widehat g_{k_1}}\widehat g_{k_2}\widehat{f-g}_{k_3},
\end{aligned}
\]
we obtain
\[
\begin{aligned}
&\|I_1^\tau(f)-I_1^\tau(g)\|_{H^r}\\
&\le \tau
\biggl[
\sum_{k\in\mathbb Z}\langle k\rangle^{2r}
\biggl(
\sum_{\substack{k_1,k_2,k_3\in\mathbb Z\\
k=-k_1+k_2+k_3\\k_1\ne k_2}}
\frac{1}{(k_2-k_1)^2}
|\widehat{f-g}_{k_1}|
|\widehat f_{k_2}|
|\widehat f_{k_3}|
\biggr)^2
\biggr]^{1/2}\\
&\quad+\tau
\biggl[
\sum_{k\in\mathbb Z}\langle k\rangle^{2r}
\biggl(
\sum_{\substack{k_1,k_2,k_3\in\mathbb Z\\
k=-k_1+k_2+k_3\\k_1\ne k_2}}
\frac{1}{(k_2-k_1)^2}
|\widehat g_{k_1}|
|\widehat{f-g}_{k_2}|
|\widehat f_{k_3}|
\biggr)^2
\biggr]^{1/2}\\
&\quad+\tau
\biggl[
\sum_{k\in\mathbb Z}\langle k\rangle^{2r}
\biggl(
\sum_{\substack{k_1,k_2,k_3\in\mathbb Z\\
k=-k_1+k_2+k_3\\k_1\ne k_2}}
\frac{1}{(k_2-k_1)^2}
|\widehat g_{k_1}|
|\widehat g_{k_2}|
|\widehat{f-g}_{k_3}|
\biggr)^2
\biggr]^{1/2}.
\end{aligned}
\]
To estimate the first term, define the auxiliary functions
\[
\xi(x)=\sum_{k\in\mathbb Z}
|\widehat{f-g}_k|e^{ikx},
\quad
\eta(x)=\sum_{k\in\mathbb Z}
|\widehat f_k|e^{ikx},
\]
which satisfy
\[
\|\xi\|_{H^r}=\|f-g\|_{H^r},
\quad
\|\eta\|_{H^r}=\|f\|_{H^r}.
\]
Then, we obtain
\begin{align*}
&\biggl[
\sum_{k\in\mathbb Z}\langle k\rangle^{2r}
\biggl(
\sum_{\substack{k_1,k_2,k_3\in\mathbb Z\\
k=-k_1+k_2+k_3\\k_1\ne k_2}}
\frac{1}{(k_2-k_1)^2}
|\widehat{f-g}_{k_1}|
|\widehat f_{k_2}|
|\widehat f_{k_3}|
\biggr)^2
\biggr]^{1/2}\\
&=\left\|
(-\Delta)^{-1}(\overline{\xi}\eta)\eta
\right\|_{H^r}\lesssim\|\xi\|_{H^r}\|\eta\|_{H^r}^2=\|f-g\|_{H^r}\|f\|_{H^r}^2.
\end{align*}
The other two terms are estimated similarly.
Combining these estimates yields
\[
\begin{aligned}
&\|I_1^\tau(f)-I_1^\tau(g)\|_{H^r}\\
&\le c\tau
\bigl(\|f\|_{H^r}^2+\|f\|_{H^r}\|g\|_{H^r}+\|g\|_{H^r}^2\bigr)
\|f-g\|_{H^r}.
\end{aligned}
\]
Finally, by the definition of $\Phi_\tau^{\mathrm{LRI1}}$
and the isometry of $e^{i\tau\frac{\Delta}{2}}$ in $H^r$, we obtain
\begin{align*}
\|\Phi_\tau^{\mathrm{LRI1}}(f)-\Phi_\tau^{\mathrm{LRI1}}(g)\|_{H^r}&\le \|f-g\|_{H^r}
+\|I_1^\tau(f)-I_1^\tau(g)\|_{H^r}\\
&\le
\bigl[1+\tau c(\|f\|_{H^r},\|g\|_{H^r})\bigr]\|f-g\|_{H^r},
\end{align*}
which completes the proof.
\end{proof}

The following estimate from \cite{LW} will be used in the local error analysis.
\begin{lemma}
\label{lem:sobolev_product}
Let $r>1/2$, $f\in H^r(\mathbb T)$, and
$g\in H^{r-1}(\mathbb T)$. Then
\[
\|fg\|_{H^{r-1}}
\lesssim \|f\|_{H^r}\|g\|_{H^{r-1}}.
\]
\end{lemma}

\begin{lemma}[Local error]
\label{lem:local_error1}
Let $r>1/2$ and assume that the exact solution
$\psi\in L^{\infty}([0,T];H^r(\mathbb{T}))$. Then
\[
\left\|
\psi(t_{n+1})
-\Phi_\tau^{\mathrm{LRI1}}\bigl(\psi(t_n)\bigr)
\right\|_{H^r}
\le c(M_r)\tau^2,
\quad 0\le t_n<t_{n+1}\le T.
\]
\end{lemma}
\begin{proof}
By \eqref{eq:expansion1} and \eqref{eq:I^tau_decomposition1}, we have
\[
\|\psi(t_{n+1})-\Phi_\tau^{\mathrm{LRI1}}(\psi(t_n))\|_{H^r}
\le c(M_r)\tau^2+\|R_1^\tau(\psi(t_n))\|_{H^r}.
\]
It therefore suffices to estimate $R_1^\tau(\psi(t_n))$.

Using $|e^{i\xi}-1|\le|\xi|$ for $\xi\in\mathbb R$, we obtain
\[
\left|
\int_0^\tau e^{isk_3(k_2-k_1)}
\bigl(e^{isk_1(k_1-k_2)}-1\bigr)\,ds
\right|
\le \frac{\tau^2}{2}|k_1|\,|k_2-k_1|.
\]
Consequently,
\[
\left|\widehat{R_1^\tau(\psi(t_n))}_k\right|
\le \frac{\tau^2}{2}
\sum_{\substack{k_1,k_2,k_3\in\mathbb{Z}\\k=-k_1+k_2+k_3\\k_1\ne k_2}}
\frac{|k_1|}{|k_2-k_1|}
|\widehat \psi_{k_1}(t_n)|\,|\widehat \psi_{k_2}(t_n)|\,|\widehat \psi_{k_3}(t_n)|.
\]

Next, we define the auxiliary function
\[
g(x)=\sum_{k\in\mathbb Z}|\widehat \psi_k(t_n)|e^{ikx},
\]
which satisfies $\|g\|_{H^r}=\|\psi(t_n)\|_{H^r}$.
By \cref{lem:sobolev_product}
and the algebra property of $H^r$, we obtain
\[
\begin{aligned}
\|R_1^\tau(\psi(t_n))\|_{H^r}
&\lesssim \tau^2
\bigl\|
\bigl[|\partial_x|^{-1}
((|\partial_x|\overline g)g)\bigr]g
\bigr\|_{H^r}\\
&\lesssim \tau^2
\|(|\partial_x|\overline g)g\|_{H^{r-1}}
\|g\|_{H^r}\\
&\lesssim \tau^2\|g\|_{H^r}^3
\leq c(M_r)\tau^2.
\end{aligned}
\]
Combining the above estimates completes the proof.
\end{proof}

Combining the stability and local error estimates above, we obtain the following global error bound.
\begin{theorem}[Global error in $H^r$]
\label{thm:global_error1}
Let $r>1/2$ and assume that the exact solution $\psi\in L^{\infty}([0,T];H^r(\mathbb{T}))$.
Let $\psi^n$ be the numerical solution computed by LRI1 \eqref{eq:LRI1}.
There exist constants $\tau_0, c>0$, depending only
on $r$, $T$, and $M_r$, such that, for all $0<\tau\le\tau_0$,
\[
\|\psi(t_n)-\psi^n\|_{H^r}
\le c\tau,\quad 0\leq n\tau\leq T.
\]
\end{theorem}
\begin{proof}
Define $e^n=\psi(t_n)-\psi^n$. We now establish a uniform bound on the numerical solution
by induction. Assume that
\[
\|\psi^j\|_{H^r}\le M_r+1,
\quad 0\le j\le m,
\]
which holds for $m=0$. Under this assumption, applying \cref{lem:stability1} and \cref{lem:local_error1} yields, for $0\le n\le m$, and $t_{n+1}\leq T$
\[
\begin{aligned}
\|e^{n+1}\|_{H^r}
&\le
\|\psi(t_{n+1})
-\Phi_\tau^{\mathrm{LRI1}}(\psi(t_n))\|_{H^r}+
\|\Phi_\tau^{\mathrm{LRI1}}(\psi(t_n))
-\Phi_\tau^{\mathrm{LRI1}}(\psi^n)\|_{H^r}\\
&\le
c(M_r)\tau^2+\left(1+\tau c(M_r)\right)
\|e^n\|_{H^r}.
\end{aligned}
\]
Since $e^0=0$, the discrete Gronwall inequality yields
\[
\|e^{m+1}\|_{H^r}\le c\tau,
\quad (m+1)\tau\le T,
\]
where $c$ depends only on $r$, $T$, and $M_r$.
Choosing $\tau_0>0$ sufficiently small, we obtain
\[
\|\psi^{m+1}\|_{H^r}
\le \|\psi(t_{m+1})\|_{H^r}+\|e^{m+1}\|_{H^r}
\le M_r+1,
\quad 0<\tau\le\tau_0.
\]
This completes the induction and establishes the global error bound.
\end{proof}

Next, we establish the first-order convergence in $L^2$ for solutions
bounded in $H^{1/2}$.

\begin{theorem}[Global error in $L^2$]
\label{thm:global_error1_L2}
Assume that $\psi\in L^\infty([0,T];H^{1/2}(\mathbb T))$
is the exact solution.
Let $\psi^n$ be the numerical solution computed by LRI1 \eqref{eq:LRI1}.
There exist constants $\tau_0,c>0$, depending only on $T$
and $M_{1/2}$, such that, for all $0<\tau\le\tau_0$,
\[
 \|\psi(t_n)-\psi^n\|_{L^2}\le c\tau,
 \quad 0\le n\tau\le T.
\]
\end{theorem}
\begin{proof}
We follow the proof of \cref{thm:global_error1} and give the additional
estimates needed in $L^2$.
For stability, we first note that the Fourier definition of $(-\Delta)^{-1}$ gives
\begin{equation}
\label{eq:L^inf_bound}
 \|(-\Delta)^{-1}F\|_{L^\infty}
 \le \sum_{k\ne0}\frac{|\widehat F_k|}{k^2}
 \lesssim \|F\|_{L^1}.
\end{equation}
Consequently, for $f_1,f_2,f_3\in L^2(\mathbb T)$,
H\"older's inequality gives
\begin{equation}
\label{eq:L2_ineq}
\begin{aligned}
 \bigl\|(-\Delta)^{-1}(\overline{f_1}f_2)f_3\bigr\|_{L^2}
 &\le
 \|(-\Delta)^{-1}(\overline{f_1}f_2)\|_{L^\infty}
 \|f_3\|_{L^2}\\
 &\lesssim
 \|f_1\|_{L^2}\|f_2\|_{L^2}\|f_3\|_{L^2}.
\end{aligned}
\end{equation}
Using \eqref{eq:L2_ineq} and following the same argument
as in the proof of \cref{lem:stability1}, we obtain
\begin{equation}
\label{eq:stability1_L2}
\|\Phi_\tau^{\mathrm{LRI1}}(f)
       -\Phi_\tau^{\mathrm{LRI1}}(g)\|_{L^2}
 \le
 \left[1+\tau c(\|f\|_{L^2},\|g\|_{L^2})\right]
 \|f-g\|_{L^2}.
\end{equation}

For the local error, the argument leading to \eqref{eq:expansion1}, combined
with \eqref{eq:L2_ineq}, yields a remainder bounded by $c(M_0)\tau^2$ in $L^2$. 
Following the proof of \cref{lem:local_error1}, we estimate the Fourier coefficients of $R_1^\tau(\psi(t_n))$ by using the inequality
\begin{equation*}
\frac{|k_1|}{|k_2-k_1|}
 \le 1+
 \frac{|k_1|^{1/2}|k_2|^{1/2}}{|k_2-k_1|},
 \quad k_1\ne k_2.
\end{equation*}
By the Sobolev embedding $H^{1/4}\hookrightarrow L^4$
and the bound $|\widehat F_k|\lesssim\|F\|_{L^1}$, we have
\[
\||\partial_x|^{-1}F\|_{L^4}
 \lesssim \||\partial_x|^{-1}F\|_{H^{1/4}}
 \lesssim \|F\|_{L^1}.
\]
Thus, with $g$ denoting the same auxiliary function as in the proof
of \cref{lem:local_error1}, H\"older's inequality and the inequalities above yield
\begin{align*}
 \|R_1^\tau(\psi(t_n))\|_{L^2}
 &\lesssim \tau^2\left[
 \left\||g|^2g\right\|_{L^2}
 +\left\|
 \left[|\partial_x|^{-1}
 \left(\left||\partial_x|^{1/2}g\right|^2\right)\right]g
 \right\|_{L^2}
 \right]\\
 &\lesssim \tau^2\left[
 \|g\|_{L^6}^3
 +\left\||\partial_x|^{1/2}g\right\|_{L^2}^2
 \|g\|_{L^4}
 \right]\\
 &\lesssim \tau^2\|g\|_{H^{1/2}}^3
 \le c(M_{1/2})\tau^2,
\end{align*}
where we use the embeddings
$H^{1/3}\hookrightarrow L^6$ and
$H^{1/4}\hookrightarrow L^4$.
Consequently,
\begin{equation}
\label{eq:local_error1_L2}
 \|\psi(t_{n+1})
       -\Phi_\tau^{\mathrm{LRI1}}(\psi(t_n))\|_{L^2}
 \le c(M_{1/2})\tau^2.
\end{equation}
Combining the stability \eqref{eq:stability1_L2} and local error \eqref{eq:local_error1_L2} estimates above
and arguing as in \cref{thm:global_error1}, we obtain the result.
\end{proof}

\section{A second-order low-regularity integrator}
In this section, we derive a second-order low-regularity integrator for the SPE \eqref{eq:1D} and present its convergence result.

\subsection{Construction of the scheme}
For the estimates used in the construction, let $r>1/2$ and assume
that the exact solution satisfies
$\psi\in L^\infty([0,T];H^{r+1}(\mathbb T))$.
For each time level $t_n$, we introduce
\begin{align*}
V_n(s)&:=(-\Delta)^{-1}\bigl(|\psi(t_n+s)|^2-C\bigr),\\
\widetilde V_n(s)&:=(-\Delta)^{-1}
\bigl(|e^{is\frac{\Delta}{2}}\psi(t_n)|^2-C\bigr).
\end{align*}

{Starting from Duhamel's principle \eqref{eq:Duhamel1}, we obtain
the following approximation for $0\leq s\leq\tau$:
\begin{equation}
\label{eq:approx2}
\psi(t_n+s)
=e^{is\frac{\Delta}{2}}\psi(t_n)
-is\,e^{is\frac{\Delta}{2}}
V_n(0)\psi(t_n)
+O_{H^r}(s^2),
\end{equation}
where the constant in $O_{H^r}(s^2)$ depends on $M_{r+1}$. This approximation and the subsequent estimates
\eqref{eq:VV1}--\eqref{eq:VV2_commutator}
are justified in \cref{lem:VV_estimates}.} Substituting \eqref{eq:approx2} into \eqref{eq:Duhamel1} gives
\begin{equation}
\label{eq:substituted-duhamel}
\begin{aligned}
\psi(t_n+\tau)
&=e^{i\tau\frac{\Delta}{2}}\psi(t_n)
-i\int_0^\tau e^{i(\tau-s)\frac{\Delta}{2}}
V_n(s)e^{is\frac{\Delta}{2}}\psi(t_n)\,ds\\
&\quad-\int_0^\tau s\,e^{i(\tau-s)\frac{\Delta}{2}}
V_n(s)e^{is\frac{\Delta}{2}}V_n(0)\psi(t_n)\,ds
+O_{H^r}(\tau^3).
\end{aligned}
\end{equation}
For the first integral, we use the estimate
\begin{equation}
\label{eq:VV1}
\|V_n(s)-\widetilde V_n(s)\|_{H^r}
\leq c(M_{r+1})s^2, \quad 0\leq s \leq \tau.
\end{equation}
For the second integral, we use
\begin{equation}
\label{eq:VV2}
\|V_n(s)-V_n(0)\|_{H^r}\leq c(M_{r+1})s,
\end{equation}
together with the commutator estimate
\begin{equation}
\label{eq:VV2_commutator}
 \left\|
 V_n(0)e^{is\frac{\Delta}{2}}V_n(0)\psi(t_n)
 -e^{is\frac{\Delta}{2}}V_n(0)^2\psi(t_n)
 \right\|_{H^r}\leq c(M_{r+1})s.
\end{equation}
Combining these estimates, we have
\begin{equation}
\label{eq:expansion2}
\psi(t_n+\tau)=e^{i\tau\frac{\Delta}{2}}\psi(t_n)-ie^{i\tau\frac{\Delta}{2}}I^{\tau}(\psi(t_n))-\frac{\tau^2}{2}e^{i\tau\frac{\Delta}{2}}V_n(0)^2\psi(t_n)+O_{H^r}(\tau^3),
\end{equation}
where $I^\tau(\psi)$ is defined as in \eqref{eq:I^tau}.

To obtain a second-order approximation of $I^\tau(\psi)$, we refine the integrand decomposition in \eqref{eq:integrand_decomposition1}. Setting
\[
 \alpha=sk_1(k_1-k_2),\quad
 \beta=sk_3(k_2-k_1),
\]
and using the identity (see \cite[Lemma 2.6]{KOS})
\begin{equation}
 e^{i(\alpha+\beta)}
 =e^{i\beta}+(e^{i\alpha}-1)
 +(e^{i\alpha}-1)(e^{i\beta}-1),
\end{equation}
we obtain the decomposition
\begin{equation}
\label{eq:I^tau_decomposition2}
I^\tau(\psi)=I_1^\tau(\psi)+I_2^\tau(\psi)+R_2^\tau(\psi),
\end{equation}
where $I_1^\tau(\psi)$ is defined as in \eqref{eq:I1},
\begin{equation}
\label{eq:I^tau2}
I_2^\tau(\psi)
=\sum_{k\in\mathbb Z}
\sum_{\substack{k_1,k_2,k_3\in\mathbb Z\\
k=-k_1+k_2+k_3\\k_1\ne k_2}}
\frac{\overline{\widehat\psi_{k_1}}
\widehat\psi_{k_2}\widehat\psi_{k_3}}
{(k_2-k_1)^2}\int_0^\tau
\left(e^{i\alpha}-1\right)ds\,e^{ikx},
\end{equation}
and
\begin{equation}
R_2^\tau(\psi)
=\sum_{k\in\mathbb Z}
\sum_{\substack{k_1,k_2,k_3\in\mathbb Z\\
k=-k_1+k_2+k_3\\k_1\ne k_2}}
\frac{\overline{\widehat\psi_{k_1}}
\widehat\psi_{k_2}\widehat\psi_{k_3}}
{(k_2-k_1)^2}\int_0^\tau
\left(e^{i\alpha}-1\right)
\left(e^{i\beta}-1\right)ds\,e^{ikx}.
\end{equation}

The integrand in $I_2^\tau(\psi)$ vanishes when $k_1=0$. For $k_1\neq 0$, using the relation
$2k_1(k_1-k_2)=(k_2-k_1)^2+k_1^2-k_2^2$,
we obtain
\[
 \int_0^\tau\left(e^{isk_1(k_1-k_2)}-1\right)ds
 =\frac{e^{i\frac{\tau}{2}\left((k_2-k_1)^2+k_1^2-k_2^2\right)}-1}
        {ik_1(k_1-k_2)}-\tau.
\]
Substituting this expression into \eqref{eq:I^tau2},
we can write $I_2^\tau(\psi)$ in physical space as
\begin{equation}
\label{eq:I^tau2_physical}
 \begin{aligned}
 I_2^\tau(\psi)
 &=-i\left\{
  e^{-i\tau\frac{\Delta}{2}}\partial_x^{-3}
  \left[
    \left(e^{-i\tau\frac{\Delta}{2}}\partial_x^{-1}\overline\psi\right)
    \left(e^{i\tau\frac{\Delta}{2}}\psi\right)
  \right]\right\}\psi\\
 &\quad +i\partial_x^{-3}    \left[(\partial_x^{-1}\overline\psi)\psi\right]\,\psi+\tau\partial_x^{-2}\left[(\overline\psi-\overline{\widehat\psi_0})\psi\right]\,\psi.
 \end{aligned}
\end{equation}

Combining \eqref{eq:expansion2},
\eqref{eq:I^tau_decomposition2}, \eqref{eq:I^tau11}, \eqref{eq:I^tau12} and \eqref{eq:I^tau2_physical}, and dropping the remainder terms, we obtain the second-order low-regularity integrator (LRI2) for the SPE \eqref{eq:1D}
\begin{equation}
\label{eq:LRI2}
\begin{aligned}
 \psi^{n+1}
 &=e^{i\tau\frac{\Delta}{2}}\psi^n
 -ie^{i\tau\frac{\Delta}{2}}\Biggl\{
       -\tau\widehat{\psi^n}_0\,
    \partial_x^{-2}(|\psi^n|^2) +i\bigl(\partial_x^{-3}(|\psi^n|^2)\bigr)
       (\partial_x^{-1}\psi^n)\\
&\quad
 -ie^{-i\tau\frac{\Delta}{2}}
  \left[
    \bigl(e^{i\tau\frac{\Delta}{2}}\partial_x^{-3}(|\psi^n|^2)\bigr)
    \bigl(e^{i\tau\frac{\Delta}{2}}\partial_x^{-1}\psi^n\bigr)
  \right]\\
 &\quad
 -i\left\{
  e^{-i\tau\frac{\Delta}{2}}\partial_x^{-3}
  \left[
    \bigl(e^{-i\tau\frac{\Delta}{2}}\partial_x^{-1}\overline{\psi^n}\bigr)
    \bigl(e^{i\tau\frac{\Delta}{2}}\psi^n\bigr)
  \right]\right\}\psi^n\\
 &\quad
 +i\partial_x^{-3}[(\partial_x^{-1}\overline{\psi^n})\psi^n]\,\psi^n+\tau\partial_x^{-2}
       \left[(\overline{\psi^n}-\overline{\widehat{\psi^n}_0})\psi^n\right]
       \,\psi^n
 \Biggr\}\\
 &\quad -\frac{\tau^2}{2}e^{i\tau\frac{\Delta}{2}}
       \left[\left(\partial_x^{-2}(|\psi^n|^2)\right)^2\psi^n\right],
 \end{aligned}
\end{equation}
where $\psi^0=\psi_0$ and $\widehat{\psi^n}_0$ denotes
the zero Fourier mode of $\psi^n$. We denote the numerical flow in \eqref{eq:LRI2} by
$\Phi_\tau^{\mathrm{LRI2}}$, i.e.,
$\psi^{n+1}=\Phi_\tau^{\mathrm{LRI2}}(\psi^n)$.
The resulting LRI2 scheme is also fully explicit. 

\subsection{Convergence analysis}
In this subsection, we establish stability and local error estimates
for the LRI2 scheme \eqref{eq:LRI2} and then derive the global error bound.
\begin{lemma}[Stability]
\label{lem:stability2}
Let $r>1/2$, $f,g\in H^r(\mathbb T)$, and $0<\tau\leq1$. We have
\[
\|\Phi_\tau^{\mathrm{LRI2}}(f)-\Phi_\tau^{\mathrm{LRI2}}(g)\|_{H^r}
\leq
\bigl(1+\tau c(\|f\|_{H^r},\|g\|_{H^r})\bigr)\|f-g\|_{H^r}.
\]
\end{lemma}
\begin{proof}
By the Fourier representations of $I_1^\tau$ and $I_2^\tau$, we have
\[
\begin{aligned}
\widehat{I_1^\tau(f)+I_2^\tau(f)}_k
&=\sum_{\substack{k_1,k_2,k_3\in\mathbb Z\\
 k=-k_1+k_2+k_3\\k_1\ne k_2}}
\frac{\overline{\widehat f_{k_1}}\widehat f_{k_2}\widehat f_{k_3}}
     {(k_2-k_1)^2}\int_0^\tau
\left(e^{isk_3(k_2-k_1)}+e^{isk_1(k_1-k_2)}-1\right)\,ds.
\end{aligned}
\]
Noting that
\[
\begin{aligned}
&\left|
\int_0^\tau
\left(e^{isk_3(k_2-k_1)}+e^{isk_1(k_1-k_2)}-1\right)\,ds
\right|\leq3\tau,
\end{aligned}
\]
and following the same argument as in the proof of \cref{lem:stability1}, we obtain
\[
\begin{aligned}
&\|I_1^\tau(f)+I_2^\tau(f)-I_1^\tau(g)-I_2^\tau(g)\|_{H^r}\\
&\leq
c\tau\bigl(\|f\|_{H^r}^2+\|f\|_{H^r}\|g\|_{H^r}
+\|g\|_{H^r}^2\bigr)\|f-g\|_{H^r}.
\end{aligned}
\]
Moreover, using
$|f|^2-|g|^2=(f-g)\overline f+g\,\overline{f-g}$ and the boundedness of $\partial_x^{-2}$, we obtain
\[
\begin{aligned}
&\left\|
\bigl(\partial_x^{-2}(|f|^2)\bigr)^2f
-\bigl(\partial_x^{-2}(|g|^2)\bigr)^2g
\right\|_{H^r}\\
&=
\left\|
\bigl(\partial_x^{-2}(|f|^2)\bigr)^2(f-g)
+\partial_x^{-2}(|f|^2-|g|^2)
\bigl[\partial_x^{-2}(|f|^2)+\partial_x^{-2}(|g|^2)\bigr]g
\right\|_{H^r}\\
&\leq
c(\|f\|_{H^r}, \|g\|_{H^r})
\|f-g\|_{H^r}.
\end{aligned}
\]
Finally, by the definition of $\Phi_\tau^{\mathrm{LRI2}}$
and the isometry of $e^{i\tau\frac{\Delta}{2}}$ in $H^r$, we obtain
\[
\begin{aligned}
&\|\Phi_\tau^{\mathrm{LRI2}}(f)-\Phi_\tau^{\mathrm{LRI2}}(g)\|_{H^r}\\
&\leq\|f-g\|_{H^r}
+\|I_1^\tau(f)+I_2^\tau(f)-I_1^\tau(g)-I_2^\tau(g)\|_{H^r}\\
&\quad+\frac{\tau^2}{2}
\left\|
\bigl(\partial_x^{-2}(|f|^2)\bigr)^2f
-\bigl(\partial_x^{-2}(|g|^2)\bigr)^2g
\right\|_{H^r}\\
&\leq
\left[1+\tau c(\|f\|_{H^r},\|g\|_{H^r})\right]\|f-g\|_{H^r},
\end{aligned}
\]
where we have used $\tau^2\leq\tau$ for $0<\tau\leq1$.
This completes the proof.
\end{proof}

To estimate the local error, we first establish the estimates
\eqref{eq:approx2} and \eqref{eq:VV1}--\eqref{eq:VV2_commutator} used in the construction of the LRI2 scheme.
\begin{lemma}
\label{lem:VV_estimates}
Let $r>1/2$ and assume that
$\psi\in L^\infty([0,T];H^{r+1}(\mathbb T))$.
Then, for $0\leq s\leq\tau$ and $t_{n+1}\leq T$,
the following estimates hold
\[
\begin{aligned}
\text{\rm(i)}\quad&
\left\|V_n(s)-V_n(0)\right\|_{H^r}
\leq c(M_{r+1})s,\\
\text{\rm(ii)}\quad&
\left\|
V_n(0)e^{is\frac{\Delta}{2}}V_n(0)\psi(t_n)
-e^{is\frac{\Delta}{2}}V_n(0)^2\psi(t_n)
\right\|_{H^r}
\leq c(M_{r+1})s,\\
\text{\rm(iii)}\quad&
\left\|
\psi(t_n+s)-e^{is\frac{\Delta}{2}}\psi(t_n)
+is\,e^{is\frac{\Delta}{2}}V_n(0)\psi(t_n)
\right\|_{H^r}
\leq c(M_{r+1})s^2,\\
\text{\rm(iv)}\quad&
\left\|V_n(s)-\widetilde V_n(s)\right\|_{H^r}
\leq c(M_{r+1})s^2.
\end{aligned}
\]
\end{lemma}
\begin{proof}
\text{\rm(i)}
Using Duhamel's formula and the algebra property of $H^r$, we have
\[
\begin{aligned}
\|\psi(t_n+s)-\psi(t_n)\|_{H^{r-1}}&\leq
\left\|\left(e^{is\frac{\Delta}{2}}-1\right)
\psi(t_n)\right\|_{H^{r-1}}
+\int_0^s
\|V_n(\sigma)\psi(t_n+\sigma)\|_{H^{r-1}}\,d\sigma\\
&\lesssim sM_{r+1}+sM_r^3
\leq c(M_{r+1})s.
\end{aligned}
\]
By the definition of $V_n$ and
\cref{lem:sobolev_product}, it follows that
\[
\begin{aligned}
\|V_n(s)-V_n(0)\|_{H^r}
&\lesssim\bigl\|
|\psi(t_n+s)|^2-|\psi(t_n)|^2
\bigr\|_{H^{r-1}}\\
&\lesssim \left(
\|\psi(t_n+s)\|_{H^r}+\|\psi(t_n)\|_{H^r}
\right)
\|\psi(t_n+s)-\psi(t_n)\|_{H^{r-1}}\\
&\leq c(M_{r+1})s.
\end{aligned}
\]

\text{\rm(ii)}
Let $f\in H^{r+2}$ and $g\in H^{r+1}$.
Differentiating $e^{i(s-\sigma)\frac{\Delta}{2}}f e^{i\sigma\frac{\Delta}{2}}g$
with respect to $\sigma$ and integrating over $[0,s]$, we obtain
\[
\begin{aligned}
e^{is\frac{\Delta}{2}}fg-f e^{is\frac{\Delta}{2}}g
&=\frac{i}{2}\int_0^s e^{i(s-\sigma)\frac{\Delta}{2}}\left(f_{xx}
e^{i\sigma\frac{\Delta}{2}}g+2f_x
e^{i\sigma\frac{\Delta}{2}}g_x\right)\,d\sigma.
\end{aligned}
\]
By the isometry of $e^{it\frac{\Delta}{2}}$ and the algebra property of $H^r$,
\begin{equation}
\label{eq:commutator_estimate}
\|e^{is\frac{\Delta}{2}}fg-f e^{is\frac{\Delta}{2}}g\|_{H^r}
\lesssim s\|f\|_{H^{r+2}}\|g\|_{H^{r+1}}.
\end{equation}
Taking $f=V_n(0)$ and $g=V_n(0)\psi(t_n)$, we obtain \text{\rm(ii)}.

\text{\rm(iii)}
By Duhamel's formula, we have
\[
\begin{aligned}
&\psi(t_n+s)-e^{is\frac{\Delta}{2}}\psi(t_n)+is\,e^{is\frac{\Delta}{2}}V_n(0)\psi(t_n)\\
&=-i\int_0^s e^{i(s-\sigma)\frac{\Delta}{2}}
\bigl(V_n(\sigma)-V_n(0)\bigr)\psi(t_n+\sigma)\,d\sigma\\
&\quad -i\int_0^s e^{i(s-\sigma)\frac{\Delta}{2}}V_n(0)
\bigl(\psi(t_n+\sigma)-e^{i\sigma\frac{\Delta}{2}}\psi(t_n)\bigr)\,d\sigma\\
&\quad-i\int_0^s e^{i(s-\sigma)\frac{\Delta}{2}}
\bigl(V_n(0)e^{i\sigma\frac{\Delta}{2}}\psi(t_n)-e^{i\sigma\frac{\Delta}{2}}V_n(0)\psi(t_n)\bigr)\,d\sigma.
\end{aligned}
\]
Applying (i) to the first integral, \eqref{eq:approx1}
to the second, and the estimate
\eqref{eq:commutator_estimate} with
$f=V_n(0)$ and $g=\psi(t_n)$ to the third, we obtain
\[
\|\psi(t_n+s)-e^{is\frac{\Delta}{2}}\psi(t_n)+is\,e^{is\frac{\Delta}{2}}V_n(0)\psi(t_n)\|_{H^r}
\leq c(M_{r+1})s^2.
\]

\text{\rm(iv)}
By (iii) and the algebra property of $H^r$, we have
\begin{align*}
|\psi(t_n+s)|^2-|e^{is\frac{\Delta}{2}}\psi(t_n)|^2=2s\operatorname{Re}\left(
-i\overline{e^{is\frac{\Delta}{2}}\psi(t_n)}
e^{is\frac{\Delta}{2}}V_n(0)\psi(t_n)
\right)+O_{H^r}(s^2),
\end{align*}
where the remainder constant depends on $M_{r+1}$.
Since $V_n(0)$ is real-valued, we have
\[
\begin{aligned}
\operatorname{Re}\left(
-i\overline{e^{is\frac{\Delta}{2}}\psi(t_n)}
V_n(0)e^{is\frac{\Delta}{2}}\psi(t_n)
\right)=
\operatorname{Re}\left(
-iV_n(0)|e^{is\frac{\Delta}{2}}\psi(t_n)|^2
\right)=0.
\end{aligned}
\]
Subtracting this vanishing term gives
\[
\begin{aligned}
&|\psi(t_n+s)|^2-|e^{is\frac{\Delta}{2}}\psi(t_n)|^2\\
&=2s\operatorname{Re}\left(
-i\overline{e^{is\frac{\Delta}{2}}\psi(t_n)}
\left(
e^{is\frac{\Delta}{2}}V_n(0)\psi(t_n)
-V_n(0)e^{is\frac{\Delta}{2}}\psi(t_n)
\right)\right)+O_{H^r}(s^2).
\end{aligned}
\]
Applying \eqref{eq:commutator_estimate} with
$f=V_n(0)$ and $g=\psi(t_n)$, we have
\[
\begin{aligned}
\|e^{is\frac{\Delta}{2}}V_n(0)\psi(t_n)
-V_n(0)e^{is\frac{\Delta}{2}}\psi(t_n)\|_{H^r}\lesssim
s\|V_n(0)\|_{H^{r+2}}\|\psi(t_n)\|_{H^{r+1}}
\leq c(M_{r+1})s.
\end{aligned}
\]
By the definitions of $V_n$ and $\widetilde V_n$, we conclude that
\[
\begin{aligned}
\|V_n(s)-\widetilde V_n(s)\|_{H^r}
\lesssim
\left\||\psi(t_n+s)|^2-|e^{is\frac{\Delta}{2}}\psi(t_n)|^2\right\|_{H^r}\leq c(M_{r+1})s^2.
\end{aligned}
\]
This completes the proof.
\end{proof}

\begin{lemma}[Local error]
\label{lem:local_error2}
Let $r>1/2$ and assume that the exact solution $\psi\in L^\infty([0,T];H^{r+1}(\mathbb T))$.
Then, for $0\leq t_n<t_{n+1}\leq T$,
\[
\|\psi(t_{n+1})-\Phi_\tau^{\mathrm{LRI2}}(\psi(t_n))\|_{H^r}
\leq c(M_{r+1})\tau^3.
\]
\end{lemma}
\begin{proof}
By \cref{lem:VV_estimates}, the remainder in \eqref{eq:expansion2} is bounded in $H^r$ by
$c(M_{r+1})\tau^3$. Combining \eqref{eq:expansion2} and \eqref{eq:I^tau_decomposition2}, we have
\[
\|\psi(t_{n+1})-\Phi_\tau^{\mathrm{LRI2}}(\psi(t_n))\|_{H^r}
\leq c(M_{r+1})\tau^3+\|R_2^\tau(\psi(t_n))\|_{H^r}.
\]
It remains to estimate $R_2^\tau(\psi(t_n))$. Using $|e^{i\xi}-1|\leq|\xi|$ for $\xi\in\mathbb R$, we obtain, for $k_1\ne k_2$,
\[
\begin{aligned}
&\frac{1}{(k_2-k_1)^2}
\left|\int_0^\tau
(e^{isk_1(k_1-k_2)}-1)(e^{isk_3(k_2-k_1)}-1)\,ds\right|\\
&\quad \leq\int_0^\tau s^2|k_1||k_3|\,ds
=\frac{\tau^3}{3}|k_1||k_3|.
\end{aligned}
\]
By the definition of $R_2^\tau(\psi(t_n))$, we have
\[
\left|\widehat{R_2^\tau(\psi(t_n))}_k\right|
\leq\frac{\tau^3}{3}
\sum_{\substack{k=-k_1+k_2+k_3\\k_1\ne k_2}}
|k_1||k_3||\widehat\psi_{k_1}(t_n)||\widehat\psi_{k_2}(t_n)||\widehat\psi_{k_3}(t_n)|.
\]
Using the auxiliary function $g$ defined in the proof of \cref{lem:local_error1}
and the algebra property of $H^r$, we obtain
\[
\begin{aligned}
\|R_2^\tau(\psi(t_n))\|_{H^r}
&\lesssim\tau^3\|(|\partial_x|\overline g)g(|\partial_x|g)\|_{H^r}
\lesssim\tau^3\|g\|_{H^r}\||\partial_x|g\|_{H^r}^2\\
&\lesssim\tau^3\|\psi(t_n)\|_{H^{r+1}}^3
\leq c(M_{r+1})\tau^3.
\end{aligned}
\]
Combining the above estimates completes the proof.
\end{proof}

The stability estimate in \cref{lem:stability2} and the local error estimate
in \cref{lem:local_error2} yield the following second-order convergence result.
\begin{theorem}[Global error in $H^r$]
\label{thm:global_error2}
Let $r>1/2$ and assume that the exact solution $\psi\in L^\infty([0,T];H^{r+1}(\mathbb T))$.
Let $\psi^n$ be the numerical solution computed by LRI2 (3.12).
There exist constants $\tau_0\in(0,1]$ and $c>0$,
depending only on $r$, $T$, and $M_{r+1}$, such that,
for all $0<\tau\leq\tau_0$,
\[
\|\psi(t_n)-\psi^n\|_{H^r}
\leq c\tau^2,\quad 0\leq n\tau\leq T.
\]
\end{theorem}
\begin{proof}
The result follows from \cref{lem:stability2} and 
\cref{lem:local_error2} by the same argument as in the proof of \cref{thm:global_error1}.
\end{proof}

For the $L^2$ error, second-order convergence can be established
under the assumption that the solution is bounded in $H^{5/4}$.
\begin{theorem}[Global error in $L^2$]
Assume that $\psi\in L^\infty([0,T];H^{5/4}(\mathbb T))$
is the exact solution.
Let $\psi^n$ be the numerical solution computed by LRI2 \eqref{eq:LRI2}.
There exist constants $\tau_0\in(0,1]$ and $c>0$,
depending only on $T$ and $M_{5/4}$, such that,
for all $0<\tau\le\tau_0$,
\[
 \|\psi(t_n)-\psi^n\|_{L^2}\le c\tau^2,
 \quad 0\le n\tau\le T.
\]
\end{theorem}
\begin{proof}
{We follow the proof of \cref{thm:global_error2} and give the additional estimates needed in $L^2$. For stability, using \eqref{eq:L^inf_bound} and H\"older's inequality, we obtain
\[
\begin{aligned}
\|\partial_x^{-2}(|f|^2)\|_{L^\infty}
&\lesssim \|f\|_{L^2}^2,\\
\|\partial_x^{-2}(|f|^2-|g|^2)\|_{L^\infty}
&\lesssim
(\|f\|_{L^2}+\|g\|_{L^2})\|f-g\|_{L^2}.
\end{aligned}
\]
Combining these estimates with the decomposition used in
the proof of \cref{lem:stability2}, we have
\[
\begin{aligned}
&\left\|
\bigl(\partial_x^{-2}(|f|^2)\bigr)^2f
-\bigl(\partial_x^{-2}(|g|^2)\bigr)^2g
\right\|_{L^2}\\
&\qquad\leq
c(\|f\|_{L^2},\|g\|_{L^2})\|f-g\|_{L^2}.
\end{aligned}
\]
Applying \eqref{eq:L2_ineq} and following the same argument
as in the proof of \cref{lem:stability2}, we obtain
\[
\|\Phi_\tau^{\mathrm{LRI2}}(f)
-\Phi_\tau^{\mathrm{LRI2}}(g)\|_{L^2}
\leq
\bigl[1+\tau c(\|f\|_{L^2},\|g\|_{L^2})\bigr]
\|f-g\|_{L^2}, \quad 0<\tau\leq1.
\]} 

For the local error, we first reconsider the estimates in \cref{lem:VV_estimates}. By the definition of the dual norm and the algebra property
of $H^1$, we obtain
\[
 |\langle ah,\varphi\rangle|
 =|\langle h,\overline a\,\varphi\rangle|
 \le \|h\|_{H^{-1}}\|\overline a\,\varphi\|_{H^1}
 \lesssim \|h\|_{H^{-1}}\|a\|_{H^1}\|\varphi\|_{H^1}.
\]
Taking the supremum over $\|\varphi\|_{H^1}=1$, we obtain
\[
 \|ah\|_{H^{-1}}
 \lesssim \|a\|_{H^1}\|h\|_{H^{-1}}.
\]
By the definition of $V_n(s)$ and
the embedding $H^1\hookrightarrow L^\infty$, we have
\[
 \|V_n(s)\|_{H^2}
 \lesssim \||\psi(t_n+s)|^2\|_{L^2}
 \lesssim M_1^2.
\]
Using
$\|(e^{is\Delta/2}-1)\psi(t_n)\|_{H^{-1}}\lesssim sM_1$
and Duhamel's formula, we obtain
\[
\begin{aligned}
 \|\psi(t_n+s)-\psi(t_n)\|_{H^{-1}}
 \lesssim sM_1
 +\int_0^s \|V_n(\sigma)\|_{L^\infty}
           \|\psi(t_n+\sigma)\|_{L^2}\,d\sigma
 \le c(M_1)s.
\end{aligned}
\]
Combining these estimates, we have
\begin{equation*}
\|V_n(s)-V_n(0)\|_{L^\infty}\le c(M_1)s.
\end{equation*}
Following the proof of \cref{lem:VV_estimates}(ii) and using
H\"older's inequality, we have, for $f\in H^2$ and $g\in H^1$,
\[
\begin{aligned}
&\left\|
e^{is\frac{\Delta}{2}}fg
-f e^{is\frac{\Delta}{2}}g
\right\|_{L^2}\\
&\lesssim
\int_0^s \left(
 \left\|f_{xx}\right\|_{L^2}
 \left\|e^{i\sigma\frac{\Delta}{2}}g\right\|_{L^\infty}
 +\left\|f_x\right\|_{L^\infty}
 \left\|\partial_x e^{i\sigma\frac{\Delta}{2}}g\right\|_{L^2}
\right)\,d\sigma\\
&\lesssim s\left\|f\right\|_{H^2}\left\|g\right\|_{H^1}.
\end{aligned}
\]
Taking $f=V_n(0)$ and $g=V_n(0)\psi(t_n)$, we obtain
\begin{align}
\label{eq:(iii)_L^2}
\left\|
 V_n(0)e^{is\frac{\Delta}{2}}
 V_n(0)\psi(t_n)
 -e^{is\frac{\Delta}{2}}
 V_n(0)^2\psi(t_n)
\right\|_{L^2}\le c(M_1)s.
\end{align}
Following the decomposition in the proof of
\cref{lem:VV_estimates}(iii), we obtain
\[
\left\|
 \psi(t_n+s)-e^{is\frac{\Delta}{2}}\psi(t_n)
 +is e^{is\frac{\Delta}{2}}
 V_n(0)\psi(t_n)
 \right\|_{L^2}\le c(M_1)s^2.
\]
Consequently, \eqref{eq:approx2} holds with a remainder bounded by
$c(M_1)s^2$ in $L^2$.
Following the proof of \cref{lem:VV_estimates}(iv) and using the preceding
$L^2$ expansion, H\"older's inequality, and the fact that
$V_n(0)$ is real-valued, we obtain
\[
\begin{aligned}
 &|\psi(t_n+s)|^2
 -\left|e^{is\frac{\Delta}{2}}\psi(t_n)\right|^2\\
 &=2s\operatorname{Re}\left(
 -i\,\overline{e^{is\frac{\Delta}{2}}\psi(t_n)}
 \left[
 e^{is\frac{\Delta}{2}}\left(V_n(0)\psi(t_n)\right)
 -V_n(0)e^{is\frac{\Delta}{2}}\psi(t_n)
 \right]\right)+O_{L^1}(s^2),
\end{aligned}
\]
where the remainder constant depends on $M_1$. Using the preceding estimate with $f=V_n(0)$ and
$g=\psi(t_n)$ and \eqref{eq:L^inf_bound}, we have
\[
 \left\|V_n(s)-\widetilde V_n(s)\right\|_{L^\infty}
 \lesssim
 \left\|
 |\psi(t_n+s)|^2
 -\left|e^{is\frac{\Delta}{2}}\psi(t_n)\right|^2
 \right\|_{L^1}\le c(M_1)s^2.
\]

Using these estimates in the derivation of \eqref{eq:expansion2}, we obtain
a remainder bounded by $c(M_1)\tau^3$ in $L^2$.
Combining \eqref{eq:expansion2} and \eqref{eq:I^tau_decomposition2}, we have
\[
 \left\|\psi(t_{n+1})
 -\Phi_\tau^{\mathrm{LRI2}}\left(\psi(t_n)\right)
 \right\|_{L^2}
 \le c(M_1)\tau^3
 +\left\|R_2^\tau\left(\psi(t_n)\right)\right\|_{L^2}.
\]
Using the same auxiliary function $g$ as in the proof of
\cref{lem:local_error2},  H\"older's inequality and the embeddings
$H^1\hookrightarrow L^\infty$ and $H^{1/4}\hookrightarrow L^4$,
we obtain
\[
\begin{aligned}
 \left\|R_2^\tau\left(\psi(t_n)\right)\right\|_{L^2}
 &\lesssim \tau^3
 \left\|
 \left(|\partial_x|\overline g\right)
 g\left(|\partial_x|g\right)
 \right\|_{L^2}\\
 &\lesssim \tau^3
 \left\|g\right\|_{L^\infty}
 \left\||\partial_x|g\right\|_{L^4}^2\\
 &\lesssim \tau^3
 \left\|g\right\|_{H^1}
 \left\||\partial_x|g\right\|_{H^{1/4}}^2\\
 &\le c(M_{5/4})\tau^3.
\end{aligned}
\]
Since $M_1\le M_{5/4}$, it follows that
\[
 \left\|\psi(t_{n+1})
 -\Phi_\tau^{\mathrm{LRI2}}\left(\psi(t_n)\right)
 \right\|_{L^2}
 \le c(M_{5/4})\tau^3.
\]
Combining the local error estimate with the stability estimate
and following the same argument as in \cref{thm:global_error2}, we complete the proof.
\end{proof}

\section{Numerical results}
In this section, we present numerical results to validate the
error estimates of the proposed schemes for the SPE \eqref{eq:1D}.
To highlight the advantages of the proposed schemes for low-regularity initial data, we compare the temporal errors of LRI1 \eqref{eq:LRI1},
LRI2 \eqref{eq:LRI2} and sLRI1 \eqref{eq:sLRI1} with classical splitting methods and exponential integrators.
All methods are discretized in space by the Fourier
pseudospectral method with a sufficiently fine grid so that
the spatial error is negligible relative to the temporal
discretization error.

Then, we introduce the numerical methods used for comparison. For each numerical approximation $\psi^n$, we define the numerical potential by $V^n=(-\Delta)^{-1}(|\psi^n|^2-C)$. The Lie--Trotter splitting method (Lie) reads
\begin{equation}
\label{eq:lie}
\psi^{n+1}=e^{i\tau\frac{\Delta}{2}}e^{-i\tau V^n}\psi^n.
\end{equation}
The Strang splitting method (Strang) is given by
\begin{equation}
\begin{aligned}
\label{eq:strang}
&V^{n,1}=(-\Delta)^{-1}\left(\left|e^{i\tau\frac{\Delta}{4}}\psi^n\right|^2-C\right),\\
&\psi^{n+1}=e^{i\tau\frac{\Delta}{4}}e^{-i\tau V^{n,1}}e^{i\tau\frac{\Delta}{4}}\psi^n.
\end{aligned}
\end{equation}
The first-order exponential integrator (EI1) reads
\begin{equation}
\label{eq:EI1}
\psi^{n+1}
=e^{i\tau\frac{\Delta}{2}}\psi^n
-i\tau e^{i\tau\frac{\Delta}{2}}V^n\psi^n,
\end{equation}
and its symmetric two-step counterpart (EI2) is given by
\begin{equation}
\label{eq:EI2}
\psi^{n+1}
=e^{i\tau\Delta}\psi^{n-1}
-2i\tau e^{i\tau\frac{\Delta}{2}}V^n\psi^n,
\quad n\geq1.
\end{equation}
The initial value is $\psi^0=\psi_0$, and the starting value $\psi^1$ for the two-step method is computed with the EI1 \eqref{eq:EI1}.

{For $r>1/2$, the standard error estimates establish
first-order convergence of Lie and EI1 in $H^r$
for solutions bounded in $H^{r+1}$, and second-order
convergence of Strang and EI2 in $H^r$
for solutions bounded in $H^{r+2}$.
In contrast, the corresponding convergence results
for LRI1 and LRI2 require only $H^r$ and $H^{r+1}$
regularity of the exact solution, respectively.}

We consider the one-dimensional SPE on
$\mathbb{T}=[0,2\pi]$ with periodic boundary conditions.
To generate the initial data, let $N$ be an even integer
and set $x_j=2\pi j/N$, $j=0,\ldots,N-1$. Take two
independent uniformly distributed random vectors in
$[0,1]^N$ and denote
\begin{equation*}
U_N=\operatorname{rand}(N,1)+i\operatorname{rand}(N,1).
\end{equation*}
{The initial data are constructed as
\begin{equation}
\label{eq:random_initial_data}
\psi_0(x):=
\frac{|\partial_{x,N}|^{-\theta}U_N}
{\left\||\partial_{x,N}|^{-\theta}U_N
\right\|_{L^\infty}},
\end{equation}
where $|\partial_{x,N}|^{-\theta}$ is the Fourier
multiplier defined by
\begin{equation*}
\left(|\partial_{x,N}|^{-\theta}\right)_k
=
\begin{cases}
|k|^{-\theta}, & k\neq0,\\
0, & k=0,
\end{cases}
\qquad k=-N/2,\ldots,N/2-1.
\end{equation*}
The parameter $\theta$ controls the decay of the
Fourier coefficients and hence the spatial smoothness
of the generated initial data. Larger values of
$\theta$ yield stronger damping of the high-frequency
components. For each test, the same random realization is used
for all methods and time step sizes.}

In the following experiments, we fix $N=2^{12}$ and
compute the solutions up to $T=1$. The reference
solutions are computed using sLRI1 \eqref{eq:sLRI1} with a very small time step $\tau_{\mathrm{ref}}=10^{-5}$
on the same spatial grid. The temporal errors are
measured in the $H^1$-norm and $L^2$-norm.

\begin{figure}[htbp]
\centering
\begin{minipage}{0.49\textwidth}
\centering
\includegraphics[width=\linewidth]{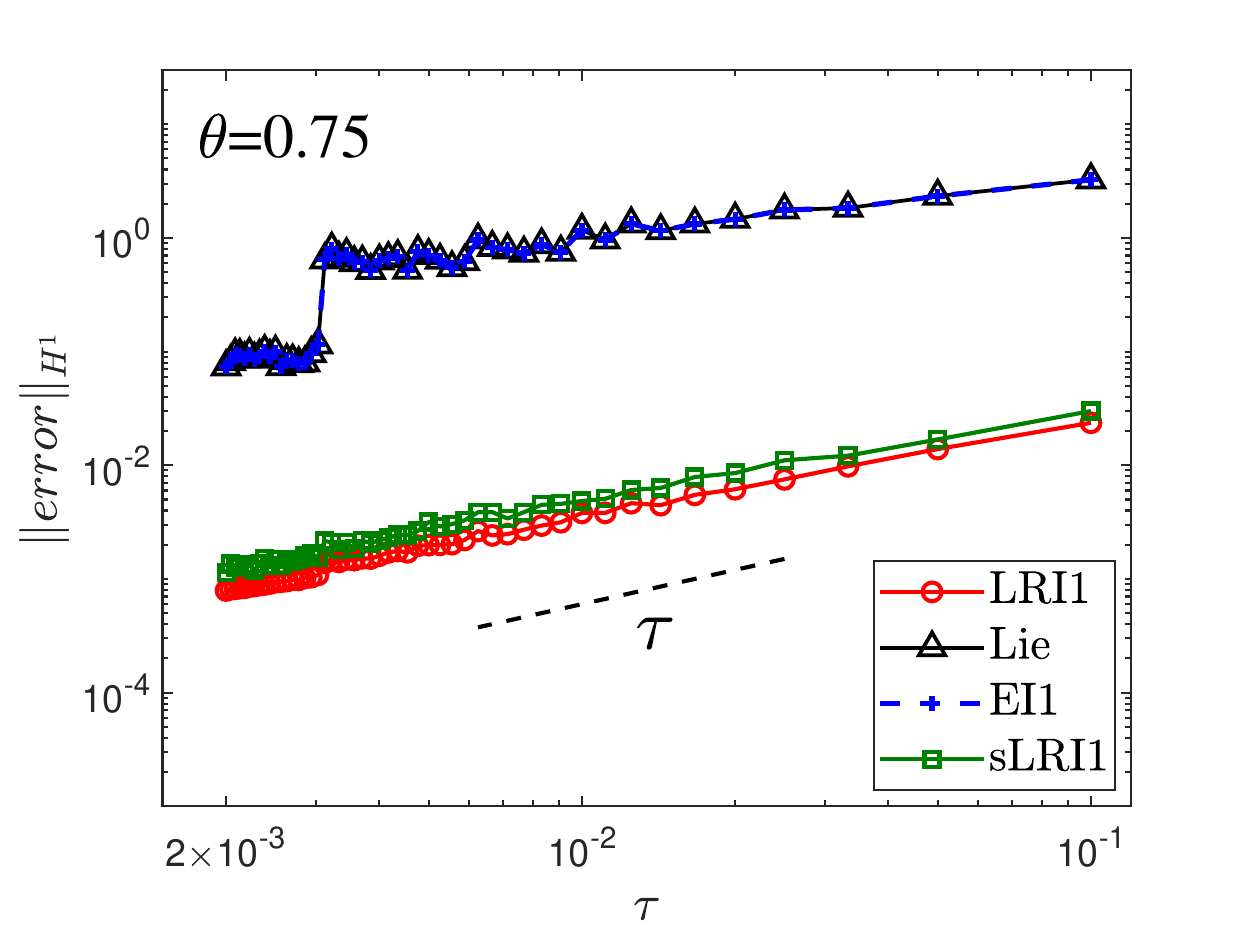}
\end{minipage}\hfill
\begin{minipage}{0.49\textwidth}
\centering
\includegraphics[width=\linewidth]{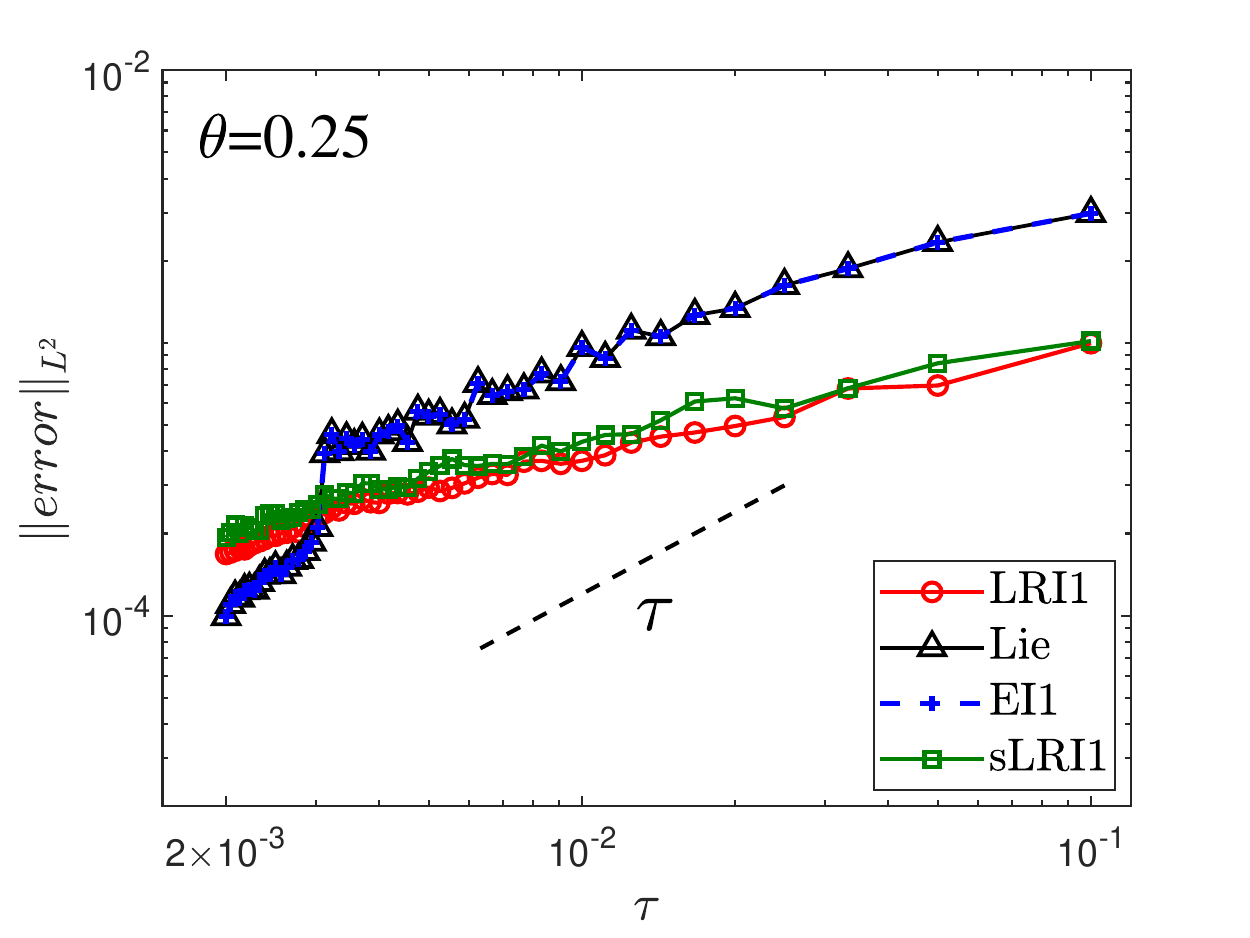}
\end{minipage}

\medskip

\begin{minipage}{0.49\textwidth}
\centering
\includegraphics[width=\linewidth]{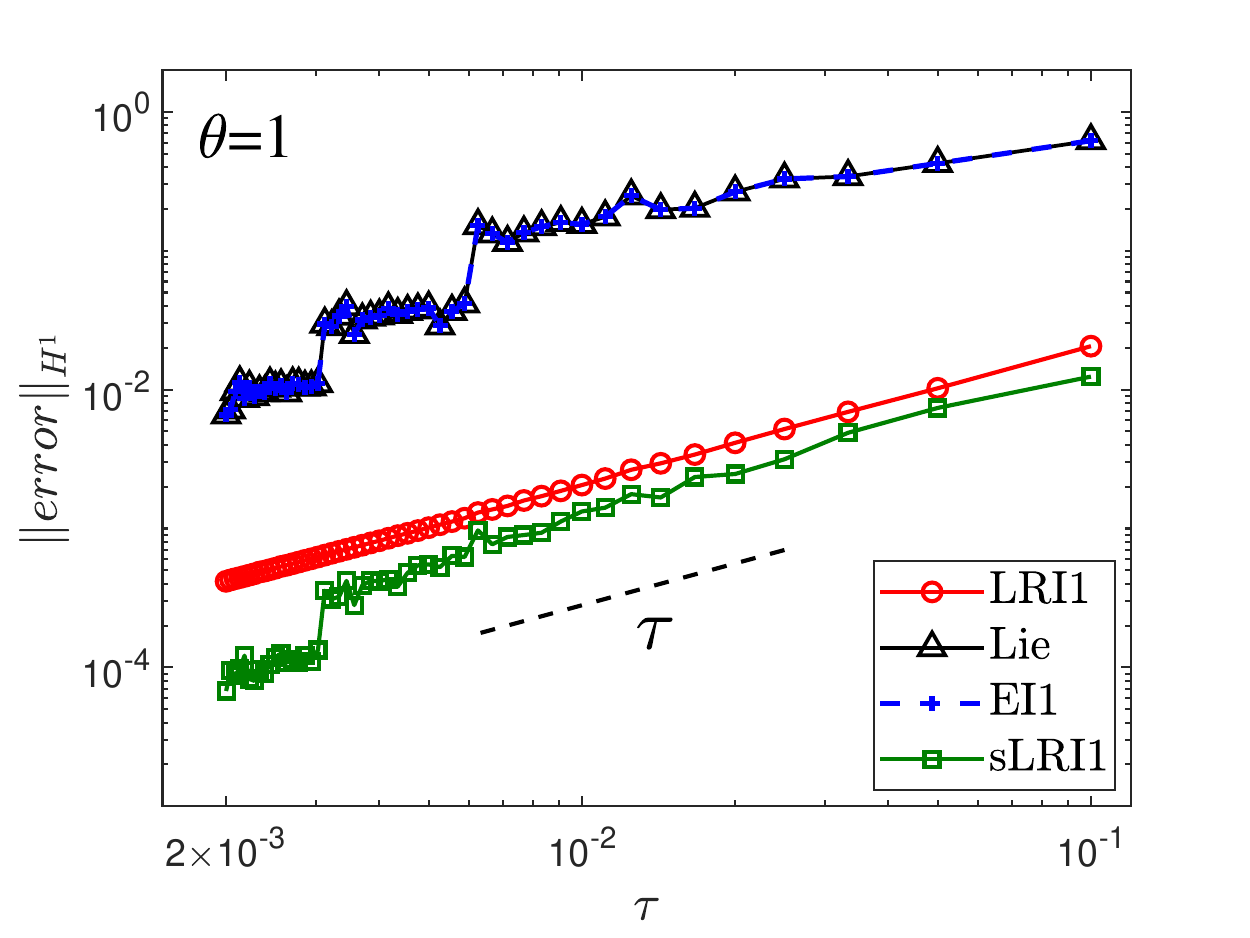}
\end{minipage}\hfill
\begin{minipage}{0.49\textwidth}
\centering
\includegraphics[width=\linewidth]{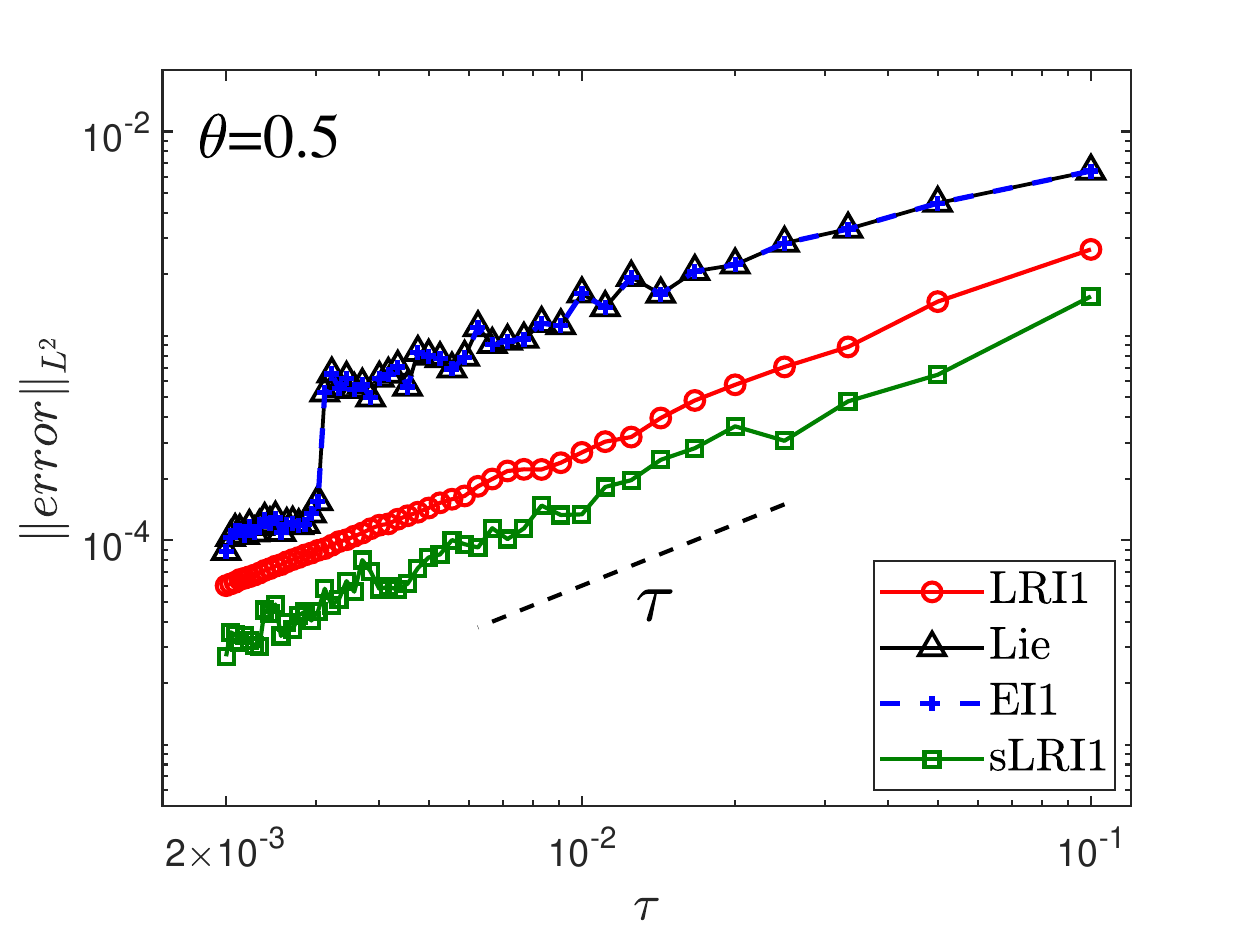}
\end{minipage}

\medskip

\begin{minipage}{0.49\textwidth}
\centering
\includegraphics[width=\linewidth]{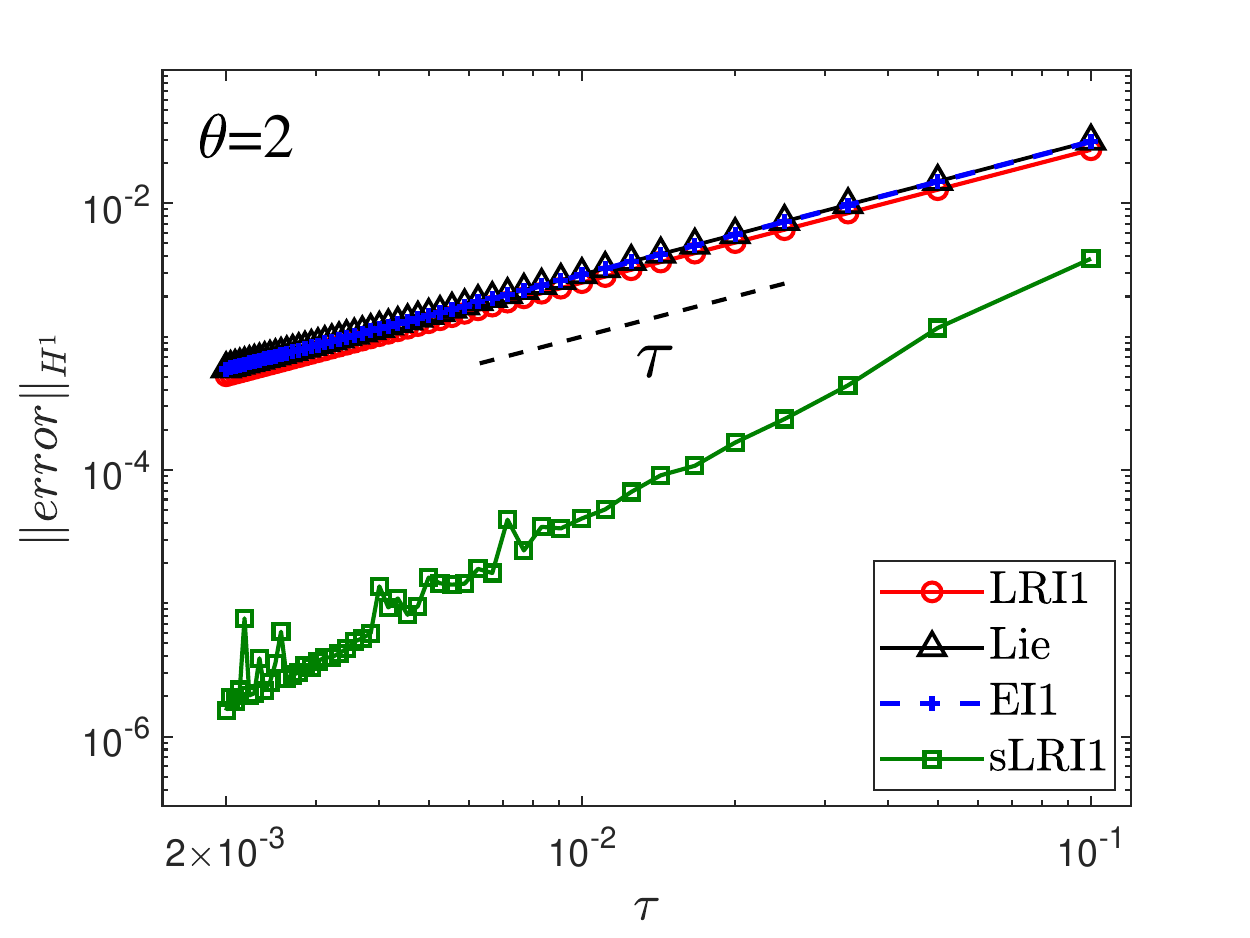}
\end{minipage}\hfill
\begin{minipage}{0.49\textwidth}
\centering
\includegraphics[width=\linewidth]{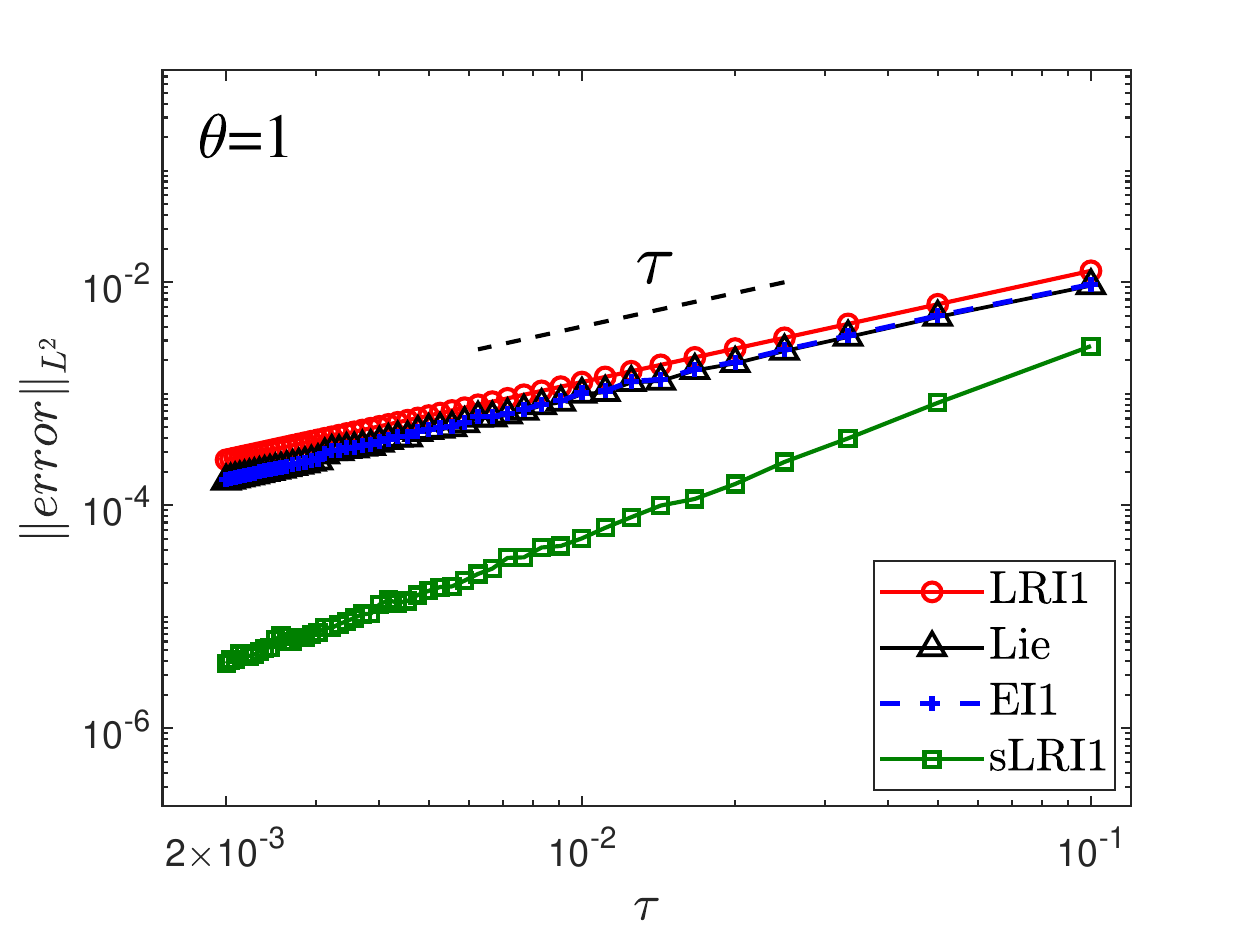}
\end{minipage}

\caption{Temporal errors in the $H^1$-norm
($\theta=0.75,1,2$, left) and the $L^2$-norm
($\theta=0.25,0.5,1$, right) of LRI1, Lie, EI1 and sLRI1
for the SPE \eqref{eq:1D} at $T=1$.}
\label{fig:first}
\end{figure}

\Cref{fig:first} presents the temporal errors of LRI1, Lie,
EI1 and sLRI1 at $T=1$ in the $H^1$-norm for
$\theta=0.75,1,2$ and in the $L^2$-norm for
$\theta=0.25,0.5,1$.
From \cref{fig:first}, we draw the following observations.

{(1)
For $\theta=0.75$ in the $H^1$-norm and $\theta=0.25$
in the $L^2$-norm, none of the four methods exhibits
a clear first-order convergence trend.
Nevertheless, for $\theta=0.75$, LRI1 and sLRI1 yield
substantially smaller $H^1$-errors than Lie and EI1.}

{(2)
For $\theta=1$ in the $H^1$-norm and $\theta=0.5$
in the $L^2$-norm, LRI1 exhibits first-order convergence,
in agreement with \cref{thm:global_error1} and
\cref{thm:global_error1_L2}, respectively, whereas
Lie and EI1 suffer from order reduction.
The errors of sLRI1 are smaller than those of LRI1,
although some oscillations are observed.}

{(3)
For $\theta=2$ in the $H^1$-norm and $\theta=1$
in the $L^2$-norm, LRI1, Lie and EI1 exhibit first-order
convergence with comparable error magnitudes.
In contrast, sLRI1 exhibits a higher observed convergence
rate and considerably smaller errors.}

{These observations support the theoretical error
estimates for LRI1 and demonstrate its ability to
retain first-order convergence under weaker regularity
conditions than Lie and EI1.
Moreover, the symmetric two-step method sLRI1
exhibits favorable error behavior, particularly for
smoother initial data.}

\begin{figure}[htbp]
\centering
\begin{minipage}{0.49\textwidth}
\centering
\includegraphics[width=\linewidth]{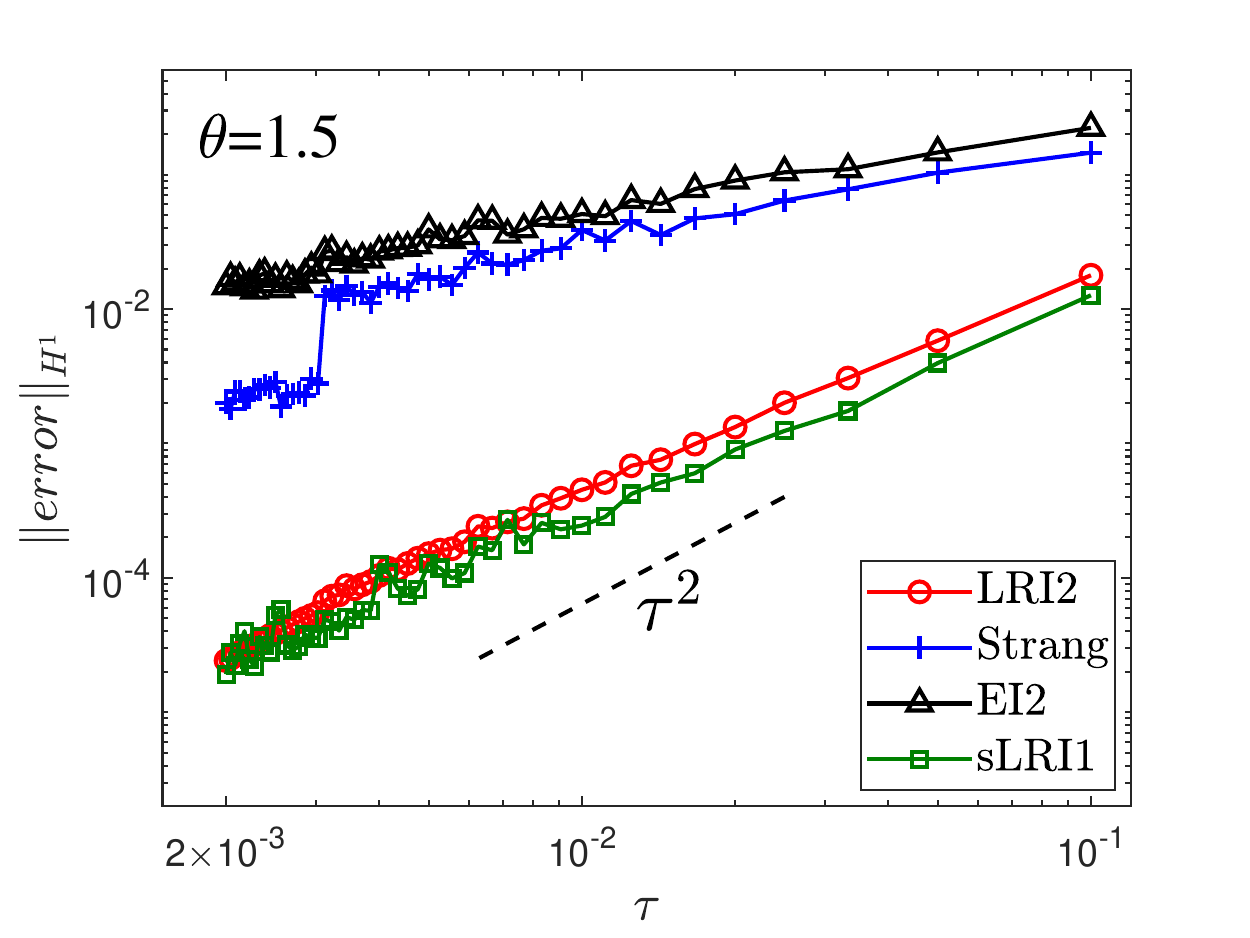}
\end{minipage}\hfill
\begin{minipage}{0.49\textwidth}
\centering
\includegraphics[width=\linewidth]{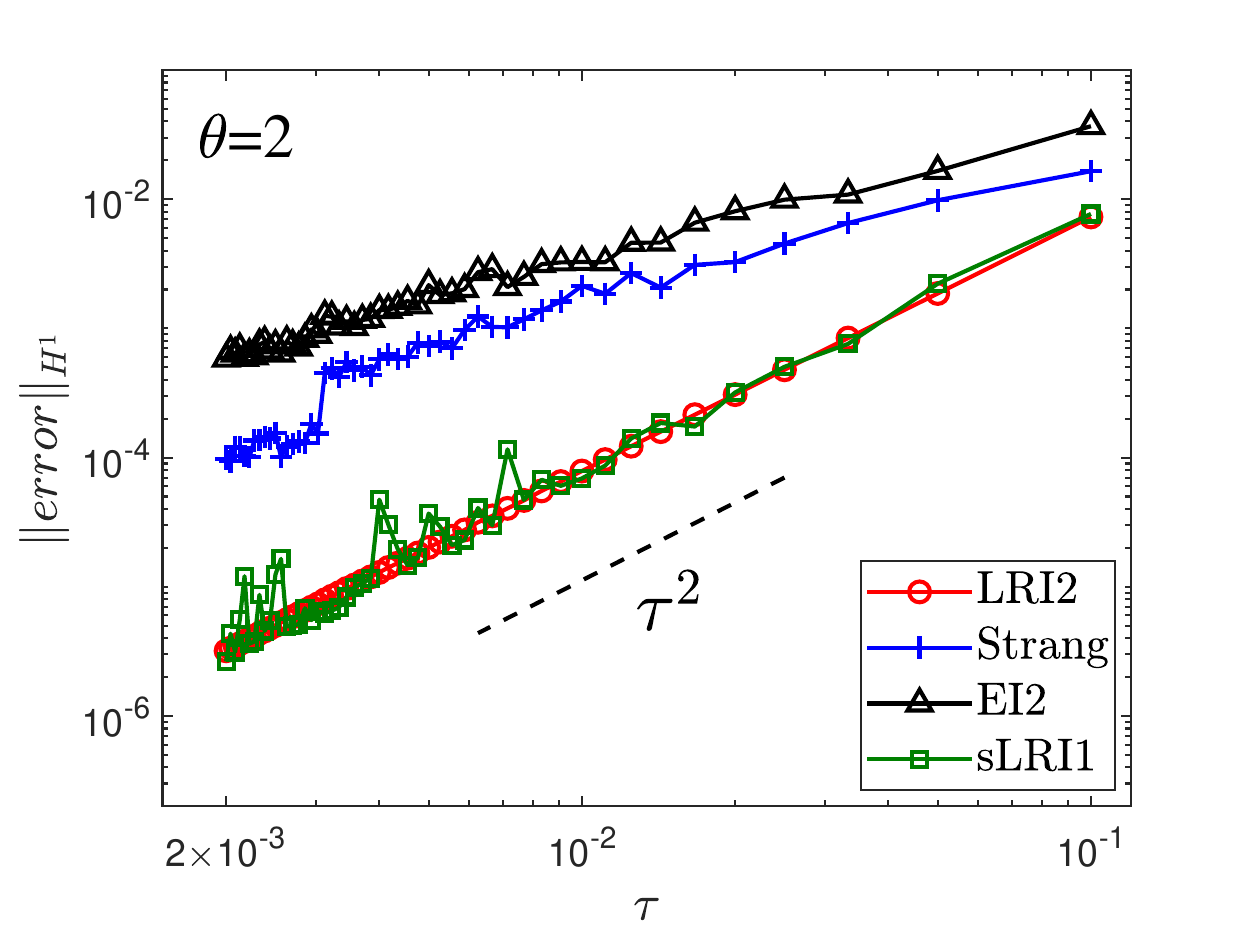}
\end{minipage}

\medskip

\begin{minipage}{0.49\textwidth}
\centering
\includegraphics[width=\linewidth]{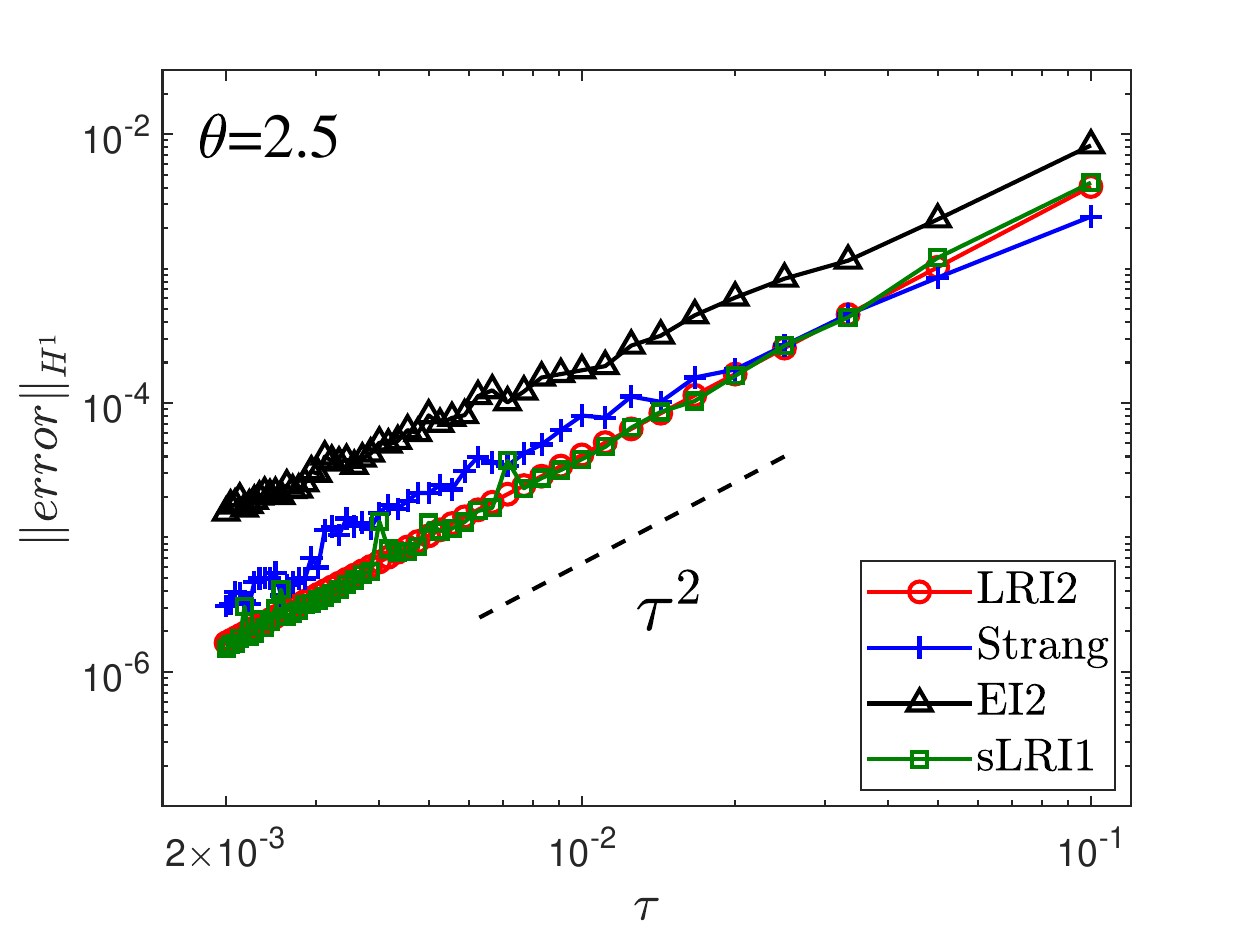}
\end{minipage}\hfill
\begin{minipage}{0.49\textwidth}
\centering
\includegraphics[width=\linewidth]{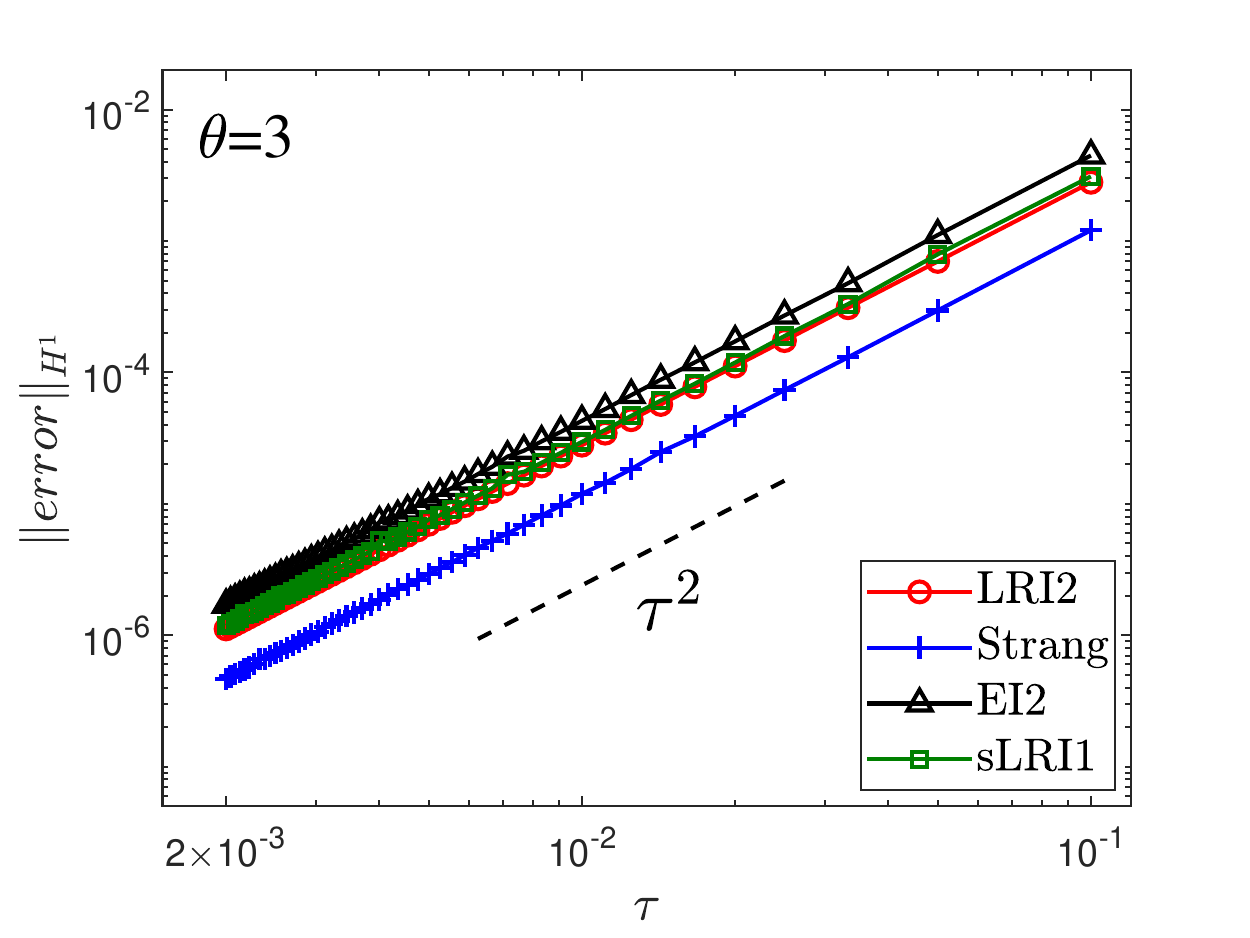}
\end{minipage}
\caption{Temporal errors in the $H^1$-norm of LRI2, Strang,
EI2 and sLRI1 for the SPE \eqref{eq:1D} at $T=1$
($\theta=1.5,2,2.5,3$).}
\label{fig:second}
\end{figure}

\Cref{fig:second} presents the temporal errors in the
$H^1$-norm of LRI2, Strang, EI2 and sLRI1 at $T=1$
for $\theta=1.5,2,2.5,3$.
We make the following observations.

{(1)
For $\theta=1.5$, none of the four methods exhibits
a clear second-order convergence trend.
Nevertheless, LRI2 and sLRI1 yield considerably
smaller errors than Strang and EI2.}

{(2)
For $\theta=2$, LRI2 exhibits second-order convergence,
in agreement with \cref{thm:global_error2}, whereas
Strang and EI2 suffer from order reduction.
The errors of sLRI1 are comparable to those of LRI2
over a substantial range of time step sizes, but
exhibit oscillations for small time step sizes.}

{(3)
For $\theta=2.5$, both LRI2 and sLRI1 exhibit
second-order convergence trends, while the convergence
rates of Strang and EI2 are lower over the range of
smaller time step sizes.}

{(4)
For $\theta=3$, all four methods exhibit second-order
convergence. In this case, Strang yields smaller errors
than the other three methods.}

{These observations support the theoretical error
estimate for LRI2 and demonstrate its advantage over
classical second-order methods for low-regularity
initial data.
In particular, LRI2 retains second-order convergence
at a lower regularity level than Strang and EI2,
while exhibiting less oscillatory error behavior
than sLRI1 in the case $\theta=2$.}

{The numerical experiments support the theoretical error estimates
and highlight the advantages of LRI1 and LRI2 for low-regularity
initial data.}

\section{Conclusion}

We developed explicit first- and second-order
low-regularity integrators for the one-dimensional periodic
Schr\"odinger--Poisson equation. The schemes are obtained by decomposing the oscillatory phase in Duhamel's formula and exploiting the smoothing effect
of the inverse Laplacian in the Poisson potential.
For $r>1/2$, the first-order method converges in $H^r$ for solutions bounded in  $H^r$, while the second-order method converges in $H^r$ for solutions bounded in $H^{r+1}$. The corresponding convergence orders in $L^2$
are established for solutions bounded in $H^{1/2}$ and
$H^{5/4}$, respectively. Numerical experiments confirm these error estimates and 
demonstrate the advantages of the proposed
integrators for low-regularity initial data, for which
classical splitting methods and exponential integrators
may suffer from order reduction.

\end{document}